\documentclass[11pt,draft]{article}
\usepackage{amsgen,amsmath,amstext,amsbsy,amsopn,amsfonts,amssymb}
\newtheorem{theorem}{Theorem}[section]
\newtheorem{lemma}[theorem]{Lemma}
\newtheorem{proposition}[theorem]{Proposition}
\newtheorem{corollary}[theorem]{Corollary}

\newtheorem{definition}{Definition}[section]
\newtheorem{note}{Note}[section]
\begin{document}
\begin{center}
\bf {\Large Stochastic evolution equations for nonlinear filtering of random fields in the presence of fractional 
Brownian sheet observation noise}
\end{center}
\begin{center}
\noindent
\rm Anna Amirdjanova\footnote{
Anna Amirdjanova, Dept of Statistics, University of Michigan, 439 West Hall, 1085 S. University Ave., Ann Arbor, MI 
48109}$^{,}$\footnote{Research supported in 
part by NSA} and Matthew Linn 
\end{center}
\begin{center}
University of Michigan 
\end{center}
{\bf Abstract.} \small 
The problem of nonlinear filtering of a random field observed in the presence 
of a noise, modeled by a persistent fractional Brownian sheet of Hurst index $(H_1,H_2)$ 
with $\displaystyle 0.5<H_1,H_2<1$, is studied and a suitable version of the 
Bayes' formula for the optimal 
filter is obtained. Two types of spatial ``fractional'' analogues of the Duncan-Mortensen-Zakai equation are also 
derived: one tracks evolution of the unnormalized optimal filter along an arbitrary ``monotone increasing'' (in the sense of 
partial ordering in $\mathbb{R}^2$) one-dimensional curve in the plane, while the other describes dynamics of the filter along paths that are 
truly two-dimensional. Although the paper deals with the two-dimensional parameter space, the presented approach and results extend to $d$-parameter random fields with arbitrary $d\geq 3$. 
\normalsize
\section{Introduction} 

An important estimation problem, arising in many engineering and physical 
systems evolving in time and space, is that of recovering a signal 
$(X_t,\;t\in\mathbb{T})$ from an observed noisy nonlinear functional 
of the signal, represented by a process $(Y_t,\;t\in\mathbb{T})$. In the  
classical mathematical filtering framework, one has $\mathbb{T}=[0,\infty)$ or 
$\mathbb{T}=[0,T]$, with $t$ interpreted as ``time", and the problem then 
is to characterize the conditional distribution of $X_t$ given the 
observation $\sigma$-field $\mathcal{F}_t^Y=\sigma\{Y_s, 0\leq s\leq t\}$, 
where the latter represents information supplied by the noisy 
observation process from time 0 up to time $t$. However, there is a 
number of interesting applications, arising, for example, in connection with 
denoising of images and video-streams,   
where the parameter space $\mathbb{T}$ has to be multidimensional, which renders   
the classical theory of nonlinear filtering inapplicable. 
The latter observation stems directly from the fact that, unlike 
$\mathbb{R}$ which permits perfect ordering, there is only partial ordering 
available in $\mathbb{R}^d$ with $d\geq 2$, thus, on the one hand, use of the multiparameter 
martingale theory in the underlying analysis is required, while, on the other hand, evolution of 
the optimal filter for $d$-parameter random fields can be studied along arbitrary $\ell$-dimensional ``monotone increasing'' 
paths with $1\leq \ell\leq d$. 

To extend the classical one-parameter nonlinear filtering theory to the multiparameter spatial filtering case, it is natural to start with the 
following observation model for a random field 
$(X_t:t\in\mathbb{T})$, with $\mathbb{T}=[0,T_1]\times\cdots\times[0,T_d]$,  
corrupted by an additive multiparameter observation noise 
$\mathcal{N}=(\mathcal{N}_t:t\in\mathbb{T})$:
\begin{equation}\label{d-model}
Y_{(t_1,\ldots,t_d)}=\int_0^{t_1}\cdots\int_0^{t_d} h\left(X_{(s_1,\ldots,s_d)}\right)
ds_1\ldots ds_d+
\mathcal{N}_{(t_1,\ldots,t_d)},\;\;\;\;(t_1,\ldots,t_d)\in\mathbb{T},  
\end{equation}  
where $h$ is a (suitably integrable) nonlinear function of the ``signal'' of interest $X$.  
Consider the 
observation $\sigma$-field  
\begin{equation}
\mathcal{F}^Y_t\equiv \mathcal{F}^Y_{(t_1,\ldots,t_d)}:=
\sigma\{Y_s:\underline{0}\prec s\prec t\},
\end{equation}
with $\underline{0}=(0,\ldots,0)\in \mathbb{T}$ and where,  
for all $t=(t_1,\ldots,t_d)$ and $s=(s_1,\ldots,s_d)$ in $\mathbb{T}$, we put 
$s\prec t$ whenever $s_i\leq t_i$ for all $i=1,\ldots,d$. Then the aim of the filtering theory is to describe the conditional distribution of the true ``signal'' of interest $X$ at ``location'' $t$, given the observation sigma-field $\mathcal{F}^Y_t$; or, equivalently, 
one can study the dynamics of $\mathbb{E}(F(X_t)\,\vert\,\mathcal{F}_t^Y)$ for a 
sufficiently rich class of test functions $F$. 

Interestingly, even in the case when observation noise $\mathcal{N}$ in~(\ref{d-model}) is a {\it standard} 
two-parameter Wiener  
sheet, assumed to be independent of $X$ ($d=2$ here), $h$ is square-integrable function and the signal is known to have a semimartingale structure, the actual derivation of 
evolution equations satisfied by the optimal filter is somewhat non-trivial, owing to the fact 
that the multiparameter martingale theory is significantly more complicated than the classical one and many ``standard'' 
martingale tools available in the one-parameter case are no longer applicable in the multiparameter setting. 
This formulation of the nonlinear spatial filtering problem (i.e. with standard Wiener sheet observation noise) 
has been studied in~\cite{bib:korezlioglu1} and \cite{bib:korezlioglu2}, where several types of stochastic partial differential equations governing the unnormalized optimal filter were obtained. However, for the case of other types of 
continuous multiparameter random fields $\mathcal{N}$ driving the observation field $Y$, no mathematical theory of optimal nonlinear filtering currently exists. 

The goal of the present paper is to study the above problem of nonlinear filtering of semimartingale random fields 
(with $d\geq 2$) but in the presence of a long-memory observation noise $\mathcal{N}$, where the latter is modelled by a persistent {\it fractional} Brownian  sheet. 

The paper is organized as follows. The remainder of Section~1 is devoted to 
two topics of interest: i) preliminaries on multiparameter martingales, which will be useful to us in Section~2; 
ii) properties of fractional Brownian sheet, plus a number of relevant results from fractional calculus. 
Section~2 presents the main theorems of the paper. Namely, as a first step, an appropriate spatial version of the ``fractional'' 
Bayes' formula is obtained. Next, a stochastic evolution equation for the unnormalized optimal filter along a one-dimensional 
monotone ``increasing'' path is derived. Finally, a stochastic evolution equation describing dynamics of the optimal 
filter along proper two-dimensional paths in the plane is also presented. Unlike what happens in the case of {\it standard} Wiener sheet observation noise, 
the latter two evolution equations cannot, strictly speaking, be interpreted as measure-valued stochastic partial differential equations due to the effects of long memory, but they certainly represent the ``fractional'' multiparameter analogues of the 
classical Duncan-Mortensen-Zakai filtering equations. Lastly some concluding remarks are given in Section~3.

For the sake of brevity and notational convenience, from now on we will restrict our attention to the case of 
two-dimensional parameter space $\mathbb{T}$, since, although analogous techniques and results can certainly be developed for the higher-dimensional cases, the latter lead to a larger number of terms in evolution equations and more cumbersome notation throughout the derivations. 
\subsection{Two-parameter martingale theory} \label{multiparametermartingale}
First let us define the usual partial ordering $\prec$ in the positive quadrant $\mathbb{R}_{+}^{2}$, along with the following 
relations and operations: for arbitrary $a=(a_{1},a_{2})$ and $b=(b_{1},b_{2})$ in 
$\mathbb{R}_{+}^{2}$, \\
$a\prec b$ if and only if $a_{1}\leq b_{1}$ and $a_{2}\leq b_{2}$;
$a\prec\prec b$ if and only if $a_{1}< b_{1}$ and $a_{2}< b_{2}$; \\
$a \curlywedge b$ if and only if $a_{1}\leq b_{1}$ and $a_{2}\geq b_{2}$;
$a\wedge b :=(\min(a_{1},b_{1}),\min(a_{2},b_{2}))$;\\
$a\vee b:=(\max(a_{1},b_{1}),\max(a_{2},b_{2}))$;
$a\odot b:=(a_{1},b_{2})$.

Given a random field $X$ with a parameter set $\mathbb{R}^2_+$, define its ``increment'' over an arbitrary 
rectangle $(z,z']:=\{x=(x_1,x_2)\in\mathbb{R}^2_+:x_1\in(z_1,z_1'],\,x_2\in (z_2,z_2']\}$, where $z \prec z'$, by  
\begin{equation}\label{X(A)}
X((z,z']):= X_{(z_{1}',z_{2}')}-X_{(z_{1},z_{2}')}-X_{(z_{1}',z_{2})}+X_{(z_{1},z_{2})}\equiv 
X_{z'}-X_{z\odot z'}-X_{z'\odot z}+X_z.
 \end{equation}
Next, for a complete probability space $(\Omega, \mathcal{F},P)$, let
$\{\mathcal{F}_{z}, z\in \mathbb{R}_{+}^{2}\}$ be a family of sub-$\sigma$-fields of $\mathcal{F}$ satisfying the following properties:\smallskip\\
$(F1)$ if $z\prec z'$ then $\mathcal{F}_{z}\subset \mathcal{F}_{z'}$;\\
$(F2)$ $\mathcal{F}_{0}$ contains all $P$-null sets of $\mathcal{F}$;\\
$(F3)$ for each $z\in\mathbb{R}^2_+$, $\mathcal{F}_{z}=\bigcap_{z \prec \prec z'}\mathcal{F}_{z'}$;\\
$(F4)$ for each $z=(z_1,z_2)\in\mathbb{R}^2_+$, $\mathcal{F}_{z}^{1}$ and $\mathcal{F}_{z}^{2}$ are conditionally independent
given $\mathcal{F}_{z}$,
where $\mathcal{F}_{z}^{1}$ and $\mathcal{F}_{z}^{2}$ are defined by
\[
\mathcal{F}_{z}^{1}:=\bigvee_{t\geq 0}\mathcal{F}_{(z_{1},t)}
= \sigma \Big\{ \bigcup_{t \in \mathbb{R}_{+}} \mathcal{F}_{(z_{1},t)} \Big\}\;\;\mathrm{and}\;\;
\mathcal{F}_{z}^{2}:=
\bigvee_{s\geq 0}\mathcal{F}_{(s,z_{2})}
= \sigma \Big\{ \bigcup_{s \in \mathbb{R}_{+}} \mathcal{F}_{(s,z_{2})} \Big\}.
\]
 Note that condition $(F4)$ is equivalent to the following condition 
$(F4')$: for all bounded random variables $X$ and all $z \in \mathbb{R}_{+}^{2}$, 
\[
\mathbb{E}\big\{X|\mathcal{F}_{z}\big\}=\mathbb{E}\big\{\mathbb{E}\{X|\mathcal{F}_{z}^{1}\big\}|
\mathcal{F}_{z}^{2}\big\}=\mathbb{E}\big\{\mathbb{E}\{X|\mathcal{F}_{z}^{2}\big\}|
\mathcal{F}_{z}^{1}\big\}\;\;\mathrm{a.s.}
\]
Moreover, for an arbitrary random field $X$ with independent increments, i.e.  
such that $X(A_{1}),\ldots, X(A_{n})$ are independent for all disjoint rectangles 
$A_{1},\ldots,A_{n}\subset\mathbb{R}^2_+$, the natural filtration generated by $X$ over rectangles 
satisfies condition $(F4)$ (see~\cite{bib:cw}), i.e. $\mathcal{F}_z^X:=\sigma\{X(A):A\prec z\}$ has property 
$(F4)$, where we say that $A\prec z$ if $x\prec z$ for all $x\in A$.
\begin{definition}
Let $(\mathcal{F}_{z})_{z\in\mathbb{R}^2_+}$ be a filtration satisfying $(F1)$-$(F4)$.
 The process $X=\{ X_{z}, z \in \mathbb{R}_{+}^{2}\}$ is called a {\it two-parameter martingale} with respect to 
$(\mathcal{F}_z)$ if: i) for each $z \in \mathbb{R}_{+}^{2}$, 
$X_{z}$ is adapted to $\mathcal{F}_{z}$ and integrable; and 
ii) for each $z \prec z'$, $\mathbb{E}\big(X_{z'}|\mathcal{F}_{z}\big)=X_{z}$ a.s. 
\end{definition}
\begin{definition}
Let $X=\{X_{z}: z\in \mathbb{R}_{+}^{2}\}$ be a process such that 
$X_{z}$ is integrable for all $z \in \mathbb{R}_{+}^{2}$ and let filtration
$(\mathcal{F}_{z})_{z\in\mathbb{R}^2_+}$ satisfy $(F1)$-$(F4)$. Then\\
(a) $X$ is called a {\it weak martingale} with respect to $(\mathcal{F}_z)$ if:

(i) $X$ is adapted to $(\mathcal{F}_{z})$, and (ii) $\mathbb{E}\big\{X((z,z'])|\mathcal{F}_{z}\big\}=0$ a.s. 
for all $z\prec\prec z'$.\smallskip\\
(b) $X$ is called an {\it $i$-martingale} ($i=1,2$) with respect to $(\mathcal{F}_z)$ if: 

(i) $X_{z}$ is $\mathcal{F}_{z}^{i}$-adapted, and 
 (ii) $\mathbb{E}\big\{X((z,z'])|\mathcal{F}_{z}^{i}\big\}=0$ a.s. for all $z\prec \prec z'$.\smallskip\\
(c) $X$ is called a {\it strong martingale} with respect to $(\mathcal{F}_z)$ if:

(i) $X$ is adapted to $(\mathcal{F}_z)$, 

(ii) $X$ vanishes on the axes (i.e. $X_{(0,z_{2})}=0$ and $X_{(z_{1},0)}=0$ a.s. for 
all $z_1,z_2\in\mathbb{R}_+$), 

and (iii) $\mathbb{E}\big\{X((z,z'])|\mathcal{F}_{z}^{1} \bigvee \mathcal{F}_{z}^{2} \big\}=0$ 
a.s. for all $z\prec\prec z'$.  
\end{definition}
Note that a martingale is both a 1- and a 2-martingale. The converse will also hold
(i.e. if $X$ is both a 1- and a 2-martingale, then $X$ is a two-parameter martingale),  
provided that 
$\{X_{(z_1,0)},\mathcal{F}_{(z_1,0)}^1,z_1\in\mathbb{R}_+\}$ and 
$\{X_{(0,z_2)},\mathcal{F}_{(0,z_2)}^2,z_2\in\mathbb{R}_+\}$ are both martingales. 
Also, clearly,
any martingale is a weak martingale and any strong martingale is a martingale. \smallskip

Let us say that a process $\{X_{z}\}$ is {\it right-continuous} if for a.e. $\omega$, 
$\lim_{\begin{subarray}{1} z' \rightarrow z \\ z \prec z' \end{subarray}}
X_{z'}(\omega)=X_{z}(\omega)$ for all $z\in\mathbb{R}^2_+$, and that it has {\it left limits} if, 
for a.e. $\omega$, $\lim_{\begin{subarray}{l} z'\rightarrow z \\ z'\prec\prec z\end{subarray}} X_{z'}(\omega)$ exists 
for all $z\in(\mathbb{R}_+\setminus\{0\})^2$. 
\begin{definition} Given filtration $(\mathcal{F}_z)$ satisfying properties $(F1)-(F4)$, 
a process $X=\{X_{z}, z\in \mathbb{R}_{+}^{2}\}$ is called an
{\it increasing process} if:  
(i) $X$ is right-continuous and adapted to $(\mathcal{F}_{z})$, 
(ii) $X_{z}=0$ a.s. on the axes, and 
(iii) $X(A) \ge 0$ for every rectangle $A\subset \mathbb{R}_{+}^{2}$.
\end{definition}
For our purposes, it will be sufficient to work with a bounded subset $\mathbb{T}=[0,T_1]\times[0,T_2]$ of $\mathbb{R}^2_+$
instead of all of $\mathbb{R}^2_+$. Let us 
fix an arbitrary $T=(T_1,T_2)\in\mathbb{R}^2_+$, and, for $p\geq 1$, define $\mathcal{M}^p(\mathbb{T})$ to be the class of 
all right-continuous martingales $M=\{M_z, z\prec T\}$ such that $M_z=0$ a.s. on the axes and 
$\mathbb{E}\vert M_z\vert^p<\infty$ for all $z\in\mathbb{T}$. Let $\mathcal{M}^p_c(\mathbb{T})$ and 
$\mathcal{M}^p_S(\mathbb{T})$ denote respectively 
the class of continuous and the class of strong martingales in $\mathcal{M}^p(\mathbb{T})$. \\

The following result highlights some of the fundamental differences between the classical one-parameter martingale theory 
and its multiparameter analogue. As shown by Cairoli and Walsh in~\cite{bib:cw} in the two-parameter case, for an arbitrary martingale 
$M\in\mathcal{M}^2(\mathbb{T})$, there exists an increasing process $A=\{A_z, z\in\mathbb{T}\}$ such that $\{M_z^2-A_z, 
z\in\mathbb{T}\}$ is a {\it weak} martingale. However, such an increasing process $A$ need not be unique even in the case 
of a strong martingale $M$. Nor can one in general guarantee the existence of an increasing process $A$ such that 
$\{M_z^2-A_z,z\in\mathbb{T}\}$ is a regular two-parameter martingale. 
Thus, we will agree to denote by $\langle M\rangle=\{\langle M\rangle_z,z\in\mathbb{T}\}$ {\it any} increasing process $A$ such that $M^2-A$ is a weak martingale.  Some refinements of the above weaker form of the Doob-Meyer decomposition are however possible in the case of 
strong martingales. Namely, if $M\in\mathcal{M}^2_S(\mathbb{T})$, then there exists a unique $\mathcal{F}_z^1$-predictable 
increasing process $[M]^{(1)}$ and there exists a unique $\mathcal{F}_z^2$-predictable increasing process $[M]^{(2)}$ such that 
$M^2_z-[M]^{(i)}_z$ is an $i$-martingale for $i=1,2$. As noted in~\cite{bib:cw}, for a strong martingale $M$, either $[M]^{(1)}$ or 
$[M]^{(2)}$ can serve as the process $\langle M\rangle$ above, but the question remains about whether the equality $[M]^{(1)}=[M]^{(2)}$ a.s. is true in general for a strong martingale $M$. In many interesting cases, the answer to the latter is in fact affirmative. 
For example, if $M\in\mathcal{M}^2_S(\mathbb{T})\cap\mathcal{M}^4_c(\mathbb{T})$, or if $M\in\mathcal{M}^2_S$ and $(\mathcal{F}_z)$ 
is a filtration generated by a standard two-parameter Wiener process (Wiener sheet), then $[M]^{(1)}=[M]^{(2)}$ a.s.\\

Although it is possible to develop the theory of stochastic integration in the plane with respect to general two-parameter martingales 
$M\in\mathcal{M}^2(\mathbb{T})$ and define corresponding stochastic integrals $\int \phi dM$ and $\iint\psi dMdM$ 
(see~\cite{bib:cw},\cite{bib:zakai1}), 
as well as the so-called mixed area integrals (as in~\cite{bib:zakai1}) of the form $\iint h d\mu dM$ and $\iint g dMd\mu$ (where 
$\mu_z,z\in\mathbb{T}$, is a continuous random function of bounded variation adapted to $(\mathcal{F}_z)$, and such that 
$\vert\mu\vert(\mathbb{T})\leq C$ a.s. for some constant $C<\infty$, where $\vert\mu\vert$ denotes the total variation measure 
corresponding to the signed measure that $\mu$ generates), but for the purposes of the present paper it will be sufficient to study 
such integrals in the special case when $M$ is a standard two-parameter Wiener process (i.e. a standard 
two-parameter Wiener sheet).  

Recall that if $W$ is a random measure in $\mathbb{R}^2_+$, which assigns to each Borel set $A$ a Gaussian random variable of mean 
zero and variance $\lambda(A)$, where $\lambda$ is the 2-dim Lebesgue measure, and which assigns independent random variables to 
disjoint sets, then the stochastic process $W=(W_z,z\in\mathbb{R}^2_+)$ defined by $W_z:=W(R_z)$, where $R_z:=(\underline{0},z]$ is the rectangle 
whose lower left-hand corner is the origin and whose upper right-hand corner is $z$, is called a two-parameter Wiener process or a 
Wiener sheet. Equivalently, one could define a two-parameter Wiener sheet $(W_z, z\in\mathbb{R}^2_+)$ 
as a continuous Gaussian random field on $\mathbb{R}^2_+$ with mean 0 and the covariance function given by:
\[
 \mathbb{E}(W_zW_{z'})=\min(z_1,z_1')\min(z_2,z_2'), \;\;\forall z,z'\in\mathbb{R}^2_+.
\]

Let $\{W_z,\mathcal{F}_z,z\in\mathbb{T}\}$ be a Wiener sheet. Let us introduce the following classes of integrands.  Let $\{\phi_z,z\in\mathbb{T}\}$ be a process such that the following conditions hold:\smallskip\\
(a) $\phi$ is a bimeasurable function of $(\omega,z)$,\\
(b) $\int_{\mathbb{T}}\mathbb{E}\phi_z^2dz<\infty$,\\
and for each $z\in\mathbb{T}$, 

either ($\mathrm{c}_0$) $\phi_z$ is $\mathcal{F}_z$-measurable, 

or ($\mathrm{c}_1$) $\phi_z$ 
is $\mathcal{F}_z^1$-measurable, 

or ($\mathrm{c}_2$) $\phi_z$ 
is $\mathcal{F}_z^2$-measurable. 
\begin{definition} \label{H0}
 For $i=0,1,2$, let $\mathcal{H}_i$ denote the space of $\phi$ satisfying (a),(b) and ($\mathrm{c}_i$).
\end{definition}
Then one can show that for $\phi\in\mathcal{H}_i,i=0,1,2,$ the stochastic integral $\int_{\mathbb{T}}\phi_zdW_z$ can be constructed (as in~\cite{bib:zakai2}). 
Moreover, if one defines the process
\[
(\phi\circ W)_z=\int_{R_z}\phi_{\zeta}dW_{\zeta}=\int_{\mathbb{T}}I(\zeta\prec z)\phi_{\zeta}dW_{\zeta},\;z\in\mathbb{T},
\]
then the process $\phi\circ W$ is a strong martingale for $\phi\in\mathcal{H}_0$, a 1-martingale for $\phi\in\mathcal{H}_1$ and 
a 2-martingale for $\phi\in\mathcal{H}_2$. Moreover, define  a process 
\[
 \xi_z=(\phi\circ W)_z(\psi\circ W)_z-\int_{R_z}\phi_{\zeta}\psi_{\zeta}d\zeta, \;z\in\mathbb{T}.
\]
Then $\xi=(\xi_z,z\in\mathbb{T})$ is a martingale with respect to $(\mathcal{F}_z)_{z\in\mathbb{T}}$ if $\phi,\psi\in\mathcal{H}_0$, a 1-martingale if $\phi,\psi\in\mathcal{H}_1$ and a 2-martingale 
if $\phi,\psi\in\mathcal{H}_2$. In all cases continuous versions of the above defined processes can be chosen.
\begin{definition} \label{hatH}
Let $\hat{\mathcal{H}}$ denote the space of functions $\psi(\omega,z,z')\equiv \psi_{z,z'}(\omega)$ on 
 $\Omega\times\mathbb{T}\times\mathbb{T}$ which satisfy the following conditions:\\
(\^a): $\psi$ is a measurable process and for all $z,z'\in\mathbb{T}$, $\psi_{z,z'}$ is 
$\mathcal{F}_{z\vee z'}$-measurable, and\\
(\^b): $\displaystyle\iint_{\mathbb{T}^2} I(z\curlywedge z')\mathbb{E}\{\psi^2_{z,z'}\}dzdz'<\infty$. 
\end{definition}
Then for arbitrary $\psi\in\hat{\mathcal{H}}$, the stochastic integrals
\[
 X_z:=\iint_{R_z\times R_z}\psi_{\zeta,\zeta'}dW_{\zeta}dW_{\zeta'},\,\,
  Y_z^1:=\iint_{R_z\times R_z}\psi_{\zeta,\zeta'}d\zeta dW_{\zeta'},\,\,
  Y_z^2:=\iint_{R_z\times R_z}\psi_{\zeta,\zeta'}dW_{\zeta} d\zeta'
 \]
are well-defined (as in~\cite{bib:zakai2}) for all $z\in\mathbb{T}$ and $X$, $Y^1$, $Y^2$ are respectively a martingale, an 
(adapted) 1-martingale and an (adapted) 2-martingale, and in all the cases the sample-continuous versions can be chosen. 
(See~\cite{bib:zakai2} for definitions of {\it adapted} 1- and 2-martingales.) Note also that the above double integrals are 
defined in such a way that only the values of the integrand on $z\curlywedge z'$ have an effect on each integral. 
\\\\
Finally, the following proposition will be useful to us later on. 
 \begin{proposition}\cite{bib:etemadi} \label{EKtheorem}
Let $\big\{X_{z}, \mathcal{F}_{z}; z \in\mathbb{T}\equiv (0,T_1]\times(0,T_2]\big\}$ be a 
strong martingale in $\mathcal{M}_S^2(\mathbb{T})$ satisfying either of the following conditions:
(i) $(\mathcal{F}_{z})$ is a filtration generated by a Brownian sheet; or 
(ii) $X\in\mathcal{M}_c^4(\mathbb{T})$. 
Then
\[\exp\biggr\{X_{z}-\frac{1}{2}\langle X\rangle_{z} \biggr\} \textrm{  is a martingale iff \;} 
\mathbb{E}\biggr[\exp\big\{X_{T}-\frac{1}{2}\langle X\rangle_{T}\big\}\biggr]=1,\]
where $T=(T_1,T_2)\in\mathbb{R}^2_+$ as before.
\end{proposition}
\subsection{Fractional calculus and properties of fractional Brownian sheet}\label{fractionalcalculus}
\begin{definition}
Let $\varphi (x) \in L_{1} (a,b)$ ($a,b\in\mathbb{R}$) and let $\alpha>0$.  The integrals

\begin{equation}\label{fracint+}
(I_{a+}^{\alpha}\varphi)(x) := \frac{1}{\Gamma(\alpha)} \int_{a}^{x} \frac{\varphi (t)}{(x-t)^{1-\alpha}}dt,
\quad x>a,
\end{equation}

\begin{equation}\label{fracint-}
(I_{b-}^{\alpha}\varphi)(x) := \frac{1}{\Gamma(\alpha)} \int_{x}^{b} \frac{\varphi (t)}{(t-x)^{1-\alpha}}dt, 
\quad x<a,
\end{equation}
are called left-sided and right-sided Riemann-Liouville fractional integrals of order~$\alpha$.
\end{definition}
Fractional integrals (\ref{fracint+}) and (\ref{fracint-}) can, in fact, be defined for functions $\varphi(x)\in L_1(a,b)$, 
existing almost everywhere. The following {\it formula for fractional integration by parts} is 
valid and will be useful (see~\cite{bib:samko}): For arbitrary $\varphi(x)\in L_p(a,b)$ and $\psi(x)\in L_q(a,b)$, where 
either $p^{-1}+q^{-1}\leq 1+\alpha$ and $(p,q)\in\mathbb{N}\setminus\{(1,1)\}$, or $p=q=1$ but $p^{-1}+q^{-1}<1+\alpha$, the 
following relation holds: 
\begin{equation} \label{intbyparts}
 \int_a^b \varphi(x)(I_{a+}^{\alpha}\psi)(x)dx=\int_a^b\psi(x)(I_{b-}^{\alpha}\varphi)(x)dx.
\end{equation}
Moreover, fractional integration has the following {\it semigroup property}: 
\begin{equation}
I_{a+}^{\alpha}I_{a+}^{\beta}\varphi=I_{a+}^{\alpha+\beta}\varphi,\;\;
I_{b-}^{\alpha}I_{b-}^{\beta}\varphi=I_{b-}^{\alpha+\beta}\varphi,\;\;\forall\alpha,\beta>0, 
\end{equation}
where the above equation holds for every point in $(a,b)$ if $\varphi\in C[a,b]$ and for almost all points in $(a,b)$ if 
$\varphi\in L_1(a,b)$.

\begin{definition}\label{fracder}
For functions $f(x)$ on interval $[a,b]\subset\mathbb{R}$, the expressions (if they exist) 

\begin{equation}\label{fracderiv+}
(\mathcal{D}_{a+}^{\alpha}f)(x):=\frac{1}{\Gamma (1-\alpha)}\frac{d}{dx}\int_{a}^{x} \frac{f(t)}{(x-t)^{\alpha}}dt,
\end{equation}

\begin{equation}\label{fracderiv-}
(\mathcal{D}_{b-}^{\alpha}f)(x):=-\frac{1}{\Gamma (1-\alpha)}\frac{d}{dx}\int_{x}^{b} \frac{f(t)}{(t-x)^{\alpha}}dt,
\end{equation}
where $0<\alpha<1$, are called respectively the left-handed and right-handed fractional Riemann-Liouville derivatives of order~$\alpha$. Moreover, for $\alpha\geq 1$, let $[\alpha]$ and $\{\alpha\}$ denote, respectively, the integral part and the ``fractional'' part 
of $\alpha$, $0\leq \{\alpha\}<1$, so that $\alpha=[\alpha]+\{\alpha\}$. Then the expressions (if they exist) 
\begin{equation}
(\mathcal{D}_{a+}^{\alpha}) f(x):=\frac{1}{\Gamma (n-\alpha)}\biggl(\frac{d}{dx}\biggr)^n
\int_{a}^{x} \frac{f(t)}{(x-t)^{\alpha-n+1}}dt,\;\;\mathrm{with\;\;}n=[\alpha]+1,
\end{equation}
\begin{equation}\label{fracderiv-}
(\mathcal{D}_{b-}^{\alpha}f)(x):=\frac{(-1)^n}{\Gamma (n-\alpha)}
\biggl(\frac{d}{dx}\biggr)^n\int_{x}^{b} \frac{f(t)}{(t-x)^{\alpha-n+1}}dt,\;\;\mathrm{with\;\;}n=[\alpha]+1,
\end{equation}
are the corresponding fractional derivatives of arbitrary order $\alpha\geq 1$. 
\end{definition}
For $\alpha<0$, we will also use the notation  
$(I_{a+}^{\alpha}\varphi)(x) := (\mathcal{D}_{a+}^{-\alpha} \varphi)(x)$ and
$(I_{b-}^{\alpha}\varphi)(x) := (\mathcal{D}_{b-}^{-\alpha} \varphi)(x).$
Also define $I_{a+}^0$ and $I_{b-}^0$ to be the identity operators: $I_{a+}^{0} \varphi=\varphi$ and $I_{b-}^{0}\varphi=\varphi$.
\begin{definition}
For $\alpha>0,$ let $I_{a+}^{\alpha}(L_{p})$ and $I_{b-}^{\alpha}(L_{q})$ be defined as the spaces of functions $f(x)$ and 
$g(x)$, respectively,  
of the form:
\begin{equation}
f=I_{a+}^{\alpha}\varphi \;\;\mathrm{for\;\;some\;\;}\varphi \in L_{p}(a,b), ~1\le p< \infty, \;\; and 
\end{equation}
\begin{equation}
g=I_{b-}^{\alpha}\psi \;\;\mathrm{for\;\;some\;\;}\psi \in L_{q}(a,b), ~1\le q< \infty.
\end{equation}
 \end{definition}
Then fractional integration and differentiation are reciprocal operations in the following sense: For $\alpha>0$ and 
arbitrary $\varphi\in L_1(a,b)$, $(\mathcal{D}_{a+}^{\alpha}I_{a+}^{\alpha}\varphi)(x)=\varphi(x)$ for almost all $x\in[a,b]$, 
while $(I_{a+}^{\alpha}\mathcal{D}_{a+}^{\alpha}f)(x)=f(x)$ is satisfied for $f\in I_{a+}^{\alpha}(L_1)$. Note, however, that there 
exist functions $f\notin I_{a+}^{\alpha}(L_1)$, whose fractional derivatives $\mathcal{D}_{a+}^{\alpha}f$ exist. 

\begin{definition}\label{Holder}
Let $\mathcal{X}$ be a finite interval.  The function $f(x)$ given on $\mathcal{X}$ is said to satisfy the 
H\"older condition of order $\lambda\in(0,1]$ on $\mathcal{X}$ (or be H\"older continuous of order $\lambda$) if there exists a constant $C>0,$ such that
\begin{equation}\label{Holdercondition}
|f(x)-f(y)| \le C |x-y|^{\lambda}
\end{equation} 
for all $x, y\in  \mathcal{X}.$
We denote by $\mathbb{H}^{\lambda}(\mathcal{X}),$ the space of functions satisfying \eqref{Holdercondition}.\\
More generally, let $\mathcal{X}$ be a finite $d$-dim rectangle. The function $f$ of $x=(x_1,\ldots,x_d)\in\mathcal{X}$ is said to satisfy 
the  H\"older condition of order $(\lambda_1,\ldots,\lambda_d)\in(0,1]^d$ on $\mathcal{X}$ (or be H\"older continuous of order 
$(\lambda_1,\ldots,\lambda_d)$) if there exists a constant $C>0,$ such that
\begin{equation}\label{Holdercondition2}
|f(x)-f(y)| \le C \left(|x_{1}-y_{1} |^{\lambda_1}+\cdots+|x_{d}-y_{d} |^{\lambda_d}\right)
\end{equation} 
for all $x,y\in\mathcal{X}$. We denote by $\mathbb{H}^{\lambda_1,\ldots,\lambda_d}(\mathcal{X}),$ the space of functions satisfying \eqref{Holdercondition2}.
\end{definition}
The following theorem will be useful to us later in the paper:
\begin{theorem}\cite{bib:samko}\label{usefulthm}
Let $f(x)=(x-a)^{-\mu}g(x)$, where $g(x) \in \mathbb{H}^{\lambda}([a,b])$, $[a,b]\subset\mathbb{R}$, $\lambda > \alpha$, 
$-\alpha< \mu < 1$.  Then,
$f(x)\in I_{a+}^{\alpha}(L_{p})$ if $\mu + \alpha < \frac{1}{p}$ for $1\le p < \infty.$
\end{theorem}

Next let us focus our attention on the properties of a {\it fractional} Brownian sheet, which will serve as a model for the 
multiparameter observation noise in the nonlinear filtering problem discussed in Section~2. The interest in studying this type of 
random field stems from the fact that it has a number of remarkable properties which make it both mathematically and practically 
interesting object, which is potentially useful in a large number of real-life applications. In fact, its one-parameter version, called  fractional Brownian motion, has recently 
become an important modelling tool in geophysical and biophysical sciences, internet traffic modelling, financial applications 
and environmental sciences. The main properties of the fractional Brownian motion (fBm) are self-similarity, 
ability to model both short and long-memory effects (depending on the value of its Hurst parameter) and its non-semimartingale and non-Markovian structure. While the latter properties often make stochastic analysis of the dynamics driven by fBm very challenging, some established close connections with the standard Wiener process through fractional 
calculus techniques (some of which were mentioned earlier in this section) provide a number of mathematical tools to make it more tractable. 

Recall that a fractional Brownian motion (fBm) with Hurst parameter $H\in(0,1)$ is a 
continuous mean zero Gaussian process ($B_{t}^{H} $, $t\in \mathbb{R}_{+}$), starting at 0 almost surely, whose  
covariance structure is given by:
\begin{equation}\label{R}
{\gamma}_{H}(s,t) := \mathbb{E}(B_{s}^{H}B_{t}^{H})=\frac{1}{2}(|s|^{2H}+|t|^{2H}-|t-s|^{2H}),  \;\; \forall
s,t\in \mathbb{R}_{+}.\end{equation} 
When $H=\frac{1}{2}$, $B_{t}^{H}$ reduces to a standard Wiener process (or standard Brownian motion). 
For $H>\frac{1}{2}$, the increments of the fBm are positively correlated and the fBm  
exhibits long range dependence (long-memory property):
$\sum_{i=1}^{\infty}\mathbb{E}\left[(B_{1}^{H}-B_0^H)(B_{i+1}^{H}-B_{i}^{H})\right]=\infty$.
When $H<\frac{1}{2}$, the increments of the fBm become negatively correlated, resulting in its short memory.
For every $H\in(0,1)$, the fBm $B^H$ is a self-similar process with self-similarity index $H$, since 
$(B^H_{ct})_{t\geq 0}\overset{d}{=}(c^HB^H_t)_{t\geq 0}$ for all constant $c>0$ (where $\overset{d}{=}$ denotes equality of two processes in 
distribution). Also, clearly, 
for any $H\ne\frac{1}{2}$, $B^H$ is not a semimartingale, implying that  
the standard techniques of stochastic calculus and stochastic integration are not directly applicable
in the fBm case. Moreover, sample paths of $B^H$ are nowhere differentiable (with probability one) but the trajectories 
are {\it H\"older continuous of any order strictly less than $H$}. Finally let us note that the fBm enjoys a number of 
fractional integral relations with respect to a standard Wiener process. For example, $B_{t}^{H}=\int_{0}^{t}K_{H}(t,s)dW_{s}$ 
for some standard Wiener process $W$, where 
\begin{equation}
K_{H}(t,s)=c_{H}\bigg(\Big(\frac{t}{s}\Big)^{H-\frac{1}{2}}(t-s)^{H-\frac{1}{2}}
-(H-\frac{1}{2})s^{\frac{1}{2}-H}\int_{s}^{t}u^{H-\frac{3}{2}}(u-s)^{H-\frac{1}{2}}du\bigg),                                         \end{equation}
where $c_{H}=\sqrt{\frac{2H\Gamma(\frac{3}{2}-H)}
{\Gamma(H+\frac{1}{2})\Gamma(2-2H)}}$.  
It is useful to note that $K_{H}(t,s)$ can also be represented by:
\begin{equation}
K_{H}(t,s)=c_H^*s^{\frac{1}{2}-H}\Big(I_{t-}^{H-\frac{1}{2}}u^{H-\frac{1}{2}}1_{[0,t]}(u) \Big)(s), 
\end{equation}
where $c_H^*=c_H\Gamma(H+\frac{1}{2})$ and 
$I_{t-}^{H-\frac{1}{2}}$ is the right-sided Riemann-Liouville fractional integral of order $H-1/2$, introduced earlier 
in this section. The above Wiener process $W$ can be reconstructed from $B^H$ via 
$W_{t}=\int_{0}^{t}K_{H}^{-1}(t,s)dB_{s}^{H}$,
where the kernel $K_{H}^{-1}$ is given by
\begin{equation}\label{K^{-1}}
K_{H}^{-1}(t,s)=c_{H}'\bigg(\Big(\frac{t}{s}\Big)^{H-\frac{1}{2}}(t-s)^{\frac{1}{2}-H}
-(H-\frac{1}{2})s^{\frac{1}{2}-H}
\int_{s}^{t}u^{H-\frac{3}{2}}(u-s)^{\frac{1}{2}-H}du\bigg),
\end{equation}
where
$
c_{H}'=\frac{1}{\Gamma(\frac{3}{2}-H)}\sqrt{\frac{\Gamma(2-2H)}
{2H\Gamma(\frac{3}{2}-H)\Gamma(H+\frac{1}{2})}}.$ The latter kernel can also be written in the following 
form:
\begin{equation}\label{K^{-1}ver2}
 K_H^{-1}(t,s)=\frac{1}{c_H^*}s^{\frac{1}{2}-H}\biggl(I_{t-}^{\frac{1}{2}-H}u^{H-\frac{1}{2}}1_{[0,t]}(u)\biggr)(s).
\end{equation}
Note that the above processes $B^H$ and 
$W$ generate the same natural filtrations.
 
The two-parameter fractional Brownian sheet (fBs) with Hurst indices $(\alpha,\beta)\in(0,1)^2$, represents the  
two-parameter analogue of fBm, and can be defined as a continuous centered Gaussian random field 
$B^{\alpha,\beta}=(B^{\alpha,\beta}_z, z\in\mathbb{R}^2_+)$, whose covariance 
structure is given by: $\mathbb{E}(B^{\alpha,\beta}_{z}B^{\alpha,\beta}_{z'})=
{\gamma}_{\alpha}(z_1,z_1'){\gamma}_{\beta}(z_2,z_2')$, with $\gamma_{\alpha}$,$\gamma_{\beta}$ defined as in~(\ref{R}) and 
$\forall z=(z_1,z_2)$,
$z'=(z_1',z_2')\in\mathbb{R}_+^2$. Naturally the fBs inherits all the remarkable properties of an fBm, while 
allowing one to introduce some new effects (like having long memory in one parameter, and short 
memory in the other, for example). In the case when both Hurst indices are greater than $1/2$, we will call 
the fBs persistent, as it displays long memory in both parameters. Note also that the fBs could also be 
equivalently defined through its integral representation with respect to a standard Wiener sheet. Namely, 
\begin{equation}\label{W_{z}}
B_{z}^{\alpha,\beta}:= \iint_{R_{z}}
K_{\alpha}(z_{1},\zeta_{1})K_{\beta}(z_{2},\zeta_{2}) dW_{(\zeta_1,\zeta_2)},\;\;z\in\mathbb{R}^2_+,
\end{equation}
where $R_{z}$ denotes the rectangle $(\underline{0},z]=(0,z_1]\times(0,z_2]$ (as before) 
and where $(W_{\zeta})_{\zeta\in\mathbb{R}^2_+}$ is a standard Wiener sheet.
On the other hand, if we let $K_{\alpha, \beta}^{-1}(z;\zeta)=
K_{\alpha}^{-1}(z_{1}, \zeta_{1})K_{\beta}^{-1}(z_{2};\zeta_{2})$, where
$K_{\alpha}^{-1}(s,t)$ is the kernel defined in \eqref{K^{-1}} (or \eqref{K^{-1}ver2}), 
then the following representation is valid: 
$W_{z}= \int_{R_{z}}
K_{\alpha,\beta}^{-1}(z;\zeta) dB_{\zeta}^{\alpha,\beta}.$ Clearly, above definition and properties 
extend to parameter spaces of dimension higher than two, but for brevity we will restrict our attention to the two-parameter 
fBs case throughout the paper. 
\section{Nonlinear filtering of random fields with persistent fractional Brownian sheet observation noise}
While in the one-parameter case the topic of optimal nonlinear filtering with fractional Gaussian observation noise has been
studied quite extensively in a variety of contexts (see e.g.,~\cite{bib:coutin},~\cite{bib:kleptsyna},~\cite{bib:mandrekar},~\cite{bib:anna},~\cite{bib:xiong},~\cite{bib:anna2}), 
there is currently no mathematical literature devoted to similar questions in the context of ``spatial'' filtering of 
multiparameter random fields. Thus, 
we expect that results presented in this paper will be of interest to  
both theoretical and applied scientists, especially in view of an increasing use of imaging technology in a number of 
fields (ranging from biomedical applications to surveillance), where the spatial structure of the underlying ``signal'' of 
interest is important and denoising and filtering of ``noisy'' images and videostreams are clearly needed. 

\subsection{Observation model with persistent fBs noise. ``Fractional-spatial'' Bayes' formula.} \label{sectionbayes} 
Consider the following observation model:
\begin{equation}\label{observation}
Y_{z}=\int_{R_{z}}g(X_{\zeta})d\zeta+B_{z}^{\alpha,\beta}, \;\;z\in\mathbb{T}\equiv [0,T_1]\times[0,T_2]
\subset \mathbb{R}_{+}^{2},
\end{equation}  
where $R_z\equiv (\underline{0},z]\equiv (0,z_1]\times(0,z_2]$, 
the signal of interest $X=(X_z,z\in\mathbb{T})$ and the observation random field $Y=(Y_z,z\in\mathbb{T})$ are measurable 
$(\mathcal{F}_z)$-adapted random fields defined on a complete filtered probability space 
$(\Omega,\mathcal{F},(\mathcal{F}_z),P)$, where the filtration $(\mathcal{F}_z)$ satisfies conditions~$(F1)$--$(F4)$ given in 
Section~1.1, 
and $B^{\alpha,\beta}=(B^{\alpha,\beta}_z,z\in\mathbb{T})$ is a fractional Brownian sheet on $(\Omega,\mathcal{F},(\mathcal{F}_z),P)$
with Hurst parameters 
$\alpha,\beta \in (\frac{1}{2},1)$ and $B^{\alpha,\beta}$ is assumed to be 
independent of the signal process $X$. Throughout the paper let us assume that the following conditions are satisfied:\\

($\mathrm{A}_1$) Function $g:\mathbb{R}\rightarrow\mathbb{R}$  
is H\"older-continuous of order $\lambda$ on any finite interval in $\mathbb{R}$, where 
$\lambda>2\max(\alpha,\beta)-1$; \\
and

($\mathrm{A}_2$) The following 
integrability condition is satisfied:

\begin{equation}\label{bassumption}
\iint_{\mathbb{T}}\left(\zeta_{1}^{\alpha-\frac{1}{2}}\zeta_{2}^{\beta-\frac{1}{2}}\right)^2\mathbb{E}\left[ 
\left(\mathcal{D}_{0+}^{\alpha-\frac{1}{2}}\otimes \mathcal{D}_{0+}^{\beta-\frac{1}{2}}
g_{\cdot}^*(X)\right)(\zeta_{1},\zeta_2)\right]^2 d\zeta_1 d\zeta_2 < \infty,
\end{equation}
where $\mathcal{D}_{0+}^{\alpha-\frac{1}{2}}$, $\mathcal{D}_{0+}^{\beta-\frac{1}{2}}$ are the fractional Riemann-Liouville derivatives defined 
in Definition~\ref{fracder}, 
\begin{equation}\label{g^*}
g_{z}^*(X)(\omega):=z_1^{\frac{1}{2}-\alpha}z_2^{\frac{1}{2}-\beta}g(X_{z}(\omega)), \;\;\forall z=(z_1,z_2)\in\mathbb{T},
\forall\omega\in\Omega,
\end{equation}
and 
``$\otimes$'' denotes the tensor product of operators.  
Namely, given a function $f:\mathbb{T}\rightarrow \mathbb{R}$ and a pair of linear operators 
$\mathcal{L}_1,\mathcal{L}_2$, defined on appropriate functions of the form $f_1:[0,T_1]\rightarrow\mathbb{R}$ and 
$f_2:[0,T_2]\rightarrow\mathbb{R}$, respectively, 
let 
\[
 \biggl(\mathcal{L}_{1}\otimes \mathcal{L}_{2}f\biggr)(z_1,z_2):=
\mathcal{L}_1(\mathcal{L}_2 f(\cdot,z_2))(z_1),\;\;\forall (z_1,z_2)\in[0,T_1]\times[0,T_2].
\]
\begin{lemma} Fix arbitrary $\alpha,\beta\in\big(\frac{1}{2},1\big)$. 
Let $h:\mathbb{T}\rightarrow\mathbb{R}$ be a H\"older continuous function of order $(\lambda_1,\lambda_2)$, 
 where $\lambda_1>\alpha-\frac{1}{2}$ and $\lambda_2>\beta-\frac{1}{2}$. Then there exists a function 
$\delta_{h}:\mathbb{T}\rightarrow\mathbb{R}$ such that $\delta_{h}\in L_2(\mathbb{T})$ and 
\begin{equation}
 \int_{[0,z_1]\times[0,z_2]} K_{\alpha,\beta}^{-1}(z_1,z_2;\zeta)h(\zeta)d\zeta=\int_{[0,z_1]\times[0,z_2]}
\delta_{h}(\zeta)d\zeta,\;\;\forall (z_1,z_2)\in\mathbb{T}. 
\end{equation}
 Moreover, if we let $h^*(z_1,z_2):=z_1^{\frac{1}{2}-\alpha}z_2^{\frac{1}{2}-\beta}h(z_1,z_2)$, then 
$\delta_{h}$ can be taken as follows:  
\begin{equation}\label{delta_h}
 \delta_{h}(z_1,z_2)=\frac{1}{c_{\alpha}^*c_{\beta}^*}\,z_1^{\alpha-\frac{1}{2}}z_2^{\beta-\frac{1}{2}}\big(\mathcal{D}_{0+}^{\alpha-\frac{1}{2}}\otimes
\mathcal{D}_{0+}^{\beta-\frac{1}{2}} h^*\big)(z_1,z_2),\;\;\forall(z_1,z_2)\in[0,T_1]\times[0,T_2].
\end{equation}
\end{lemma}
{\bf Proof}: Since linear combinations of tensor products of functions of single variable are dense in the space of 
functions of two variables, it suffices to consider the case of functions $h$ of the form $h(z_1,z_2)=
h_1(z_1)h_2(z_2)$, where $h_1\in\mathbb{H}^{\lambda_1}([0,T_1])$, $h_2\in\mathbb{H}^{\lambda_2}([0,T_2])$. By 
\eqref{K^{-1}ver2},
\begin{eqnarray*}
 \int_{[0,z_1]\times[0,z_2]}K_{\alpha,\beta}^{-1}(z_1,z_2;\zeta)h(\zeta)d\zeta&=&
\frac{1}{c_{\alpha}^*c_{\beta}^*}\int_0^{z_1}\zeta_1^{\frac{1}{2}-\alpha}\big(\mathcal{D}^{\alpha-\frac{1}{2}}_{z_1-} 
u^{\alpha-\frac{1}{2}}1_{[0,z_1]}(u)\big)(\zeta_1)h_1(\zeta_1)d\zeta_1
\\
& &\times\int_0^{z_2}\zeta_2^{\frac{1}{2}-\beta}\big(\mathcal{D}^{\beta-\frac{1}{2}}_{z_2-} 
u^{\beta-\frac{1}{2}}1_{[0,z_2]}(u)\big)(\zeta_2)h_2(\zeta_2)d\zeta_2
\end{eqnarray*}
\begin{eqnarray*}
 &=&\frac{1}{c_{\alpha}^*c_{\beta}^*}\biggl[\int_0^{z_1}\zeta_1^{\alpha-\frac{1}{2}}\big(\mathcal{D}^{\alpha-\frac{1}{2}}_{0+} 
u^{\frac{1}{2}-\alpha}h_1(u)\big)(\zeta_1)d\zeta_1\biggr]
\biggl[\int_0^{z_2}\zeta_2^{\beta-\frac{1}{2}}\big(\mathcal{D}^{\beta-\frac{1}{2}}_{0+} 
u^{\frac{1}{2}-\beta}h_2(u)\big)(\zeta_2)d\zeta_2\biggr]\\
&=&\int_{[0,z_1]\times[0,z_2]}\delta_h(\zeta_1,\zeta_2)d\zeta_1d\zeta_2,
\end{eqnarray*}
where $\delta_h$ is defined by \eqref{delta_h} and where we used the fractional differentiation by parts formula 
(see corollaries from~\eqref{intbyparts} in~\cite{bib:samko}), together with Theorem~\ref{usefulthm}. Note also that if  $h_1\in\mathbb{H}^{\lambda_1}([0,T_1])$, 
where $\lambda_1>\alpha-\frac{1}{2}$, then, by Theorem~\ref{usefulthm},   
$h_1\in I_{0+}^{\alpha-\frac{1}{2}}(L_2([0,T_1]))$, 
which implies that $\int_0^{\cdot} h_1(s)ds\in I_{0+}^{\alpha+\frac{1}{2}}(L_2([0,T_1]))$. 
If we let 
\begin{equation}
\delta_{h_1}(s):=\frac{1}{c_{\alpha}^*}s^{\alpha-\frac{1}{2}}
\big(\mathcal{D}_{0+}^{\alpha-\frac{1}{2}}v^{\frac{1}{2}-\alpha}h_1(v)\big)(s),\;\;
\forall s\in[0,T_1],
\end{equation}
then one can easily check that 
\[
 \int_0^t K_{\alpha}(t,s)\delta_{h_1}(s)ds=\int_0^t h_1(s)ds,\;\;\forall t\in[0,T_1].
\]
Since the integral operator $\mathcal{K}_{\alpha}$ associated with the kernel 
$K_{\alpha}$, 
i.e. \[
\big[\mathcal{K}_{\alpha}f\big](t)=\int_0^t K_{\alpha}(t,s)f(s)ds,
\;\;f\in L_2([0,T_1]),\] is an 
isomorphism from $L_2([0,T_1])$ onto $I_{0+}^{\alpha+\frac{1}{2}}(L_2([0,T_1]))$, 
then $f\in L_2([0,T_1])$ if and only if 
$\mathcal{K}_{\alpha}f\in I_{0+}^{\alpha+\frac{1}{2}}(L_2([0,T_1]))$. Therefore, 
$\delta_{h_1}\in L_2([0,T_1])$. Similarly one shows that if 
$h_2\in\mathbb{H}^{\lambda_2}([0,T_2])$, where $\lambda_2>\beta-\frac{1}{2}$, then 
$\delta_{h_2}\in L_2([0,T_2])$, where 
\[
\delta_{h_2}(s):=\frac{1}{c_{\beta}^*}s^{\beta-\frac{1}{2}}
\big(\mathcal{D}_{0+}^{\beta-\frac{1}{2}}v^{\frac{1}{2}-\beta}h_2(v)\big)(s),\;\;
\forall s\in[0,T_1]. 
\]
Therefore, $\delta_{h_1\otimes h_2}=\delta_{h_1}\otimes\delta_{h_2}
\in L_2([0,T_1]\times[0,T_2])$, and the required result follows. $\Box$
\begin{corollary}\label{Cor1}
 Fix $\lambda_0\in\big(\frac{\max(\alpha,\beta)-\frac{1}{2}}{\lambda},\frac{1}{2}\big)$. Suppose the signal 
$X=(X_z,z\in\mathbb{T})$ has almost surely H\"older-continuous sample paths of order $(\lambda_0,\lambda_0)$ and 
$g$ satisfies condition~($\mathrm{A}_1$). Then for almost all $\omega\in\Omega$ one can define function 
$(\delta_z(X),z\in\mathbb{T})$ by
\begin{equation}\label{delta}
 \delta_{z}(X)(\omega):=\frac{1}{c_{\alpha}^*c_{\beta}^*}
z_1^{\alpha-\frac{1}{2}}z_2^{\beta-\frac{1}{2}}\big(
\mathcal{D}_{0+}^{\alpha-\frac{1}{2}}\otimes\mathcal{D}_{0+}^{\beta-\frac{1}{2}}g^*_{\cdot}(X)(\omega)\big)(z),\;\;
z=(z_1,z_2)\in\mathbb{T},
\end{equation}
with $g^*_{\cdot}(X)$ defined by~\eqref{g^*}. Then, assuming also that~($\mathrm{A}_2$) holds,  
$\delta(X)=(\delta_{z}(X),z\in\mathbb{T})$ has the following properties: 

(i) $\delta_{\cdot}(X)(\omega)\in L_2(\mathbb{T})$ for almost all $\omega\in\Omega$ and 
$\mathbb{E}\int_{\mathbb{T}}(\delta_z(X))^2dz<\infty$;

(ii) For every rectangle $R_z=[\underline{0},z]\subset\mathbb{T}$,
\begin{equation}
 \int_{R_z} K_{\alpha,\beta}^{-1}(z;\zeta)g(X_{\zeta})d\zeta=
\int_{R_z}\delta_{\zeta}(X)d\zeta\;\;\mathrm{a.s.}
\end{equation}
\end{corollary}

From now on suppose that the assumptions of Corollary~\ref{Cor1} are satisfied. Let us introduce processes 
\begin{equation}
 W_z^Y:=\int_{R_z} K_{\alpha,\beta}^{-1}(z;\zeta)dY_{\zeta}\;\;
\mathrm{and}\;\;
W_z^B:=\int_{R_z}K_{\alpha,\beta}^{-1}(z;\zeta)dB_{\zeta}^{\alpha,\beta},\;\;z\in\mathbb{T}.
\end{equation}
Then it is easy to see that $W_z^Y=\int_{R_z}\delta_{\zeta}(X)d\zeta+W_z^B$. Next let us 
define a process $V=(V_z,z\in\mathbb{T})$ by:
\begin{equation}\label{V1}
 V_z=\exp\biggl\{-\int_{R_z}\delta_{\zeta}(X)dW_{\zeta}^B-\frac{1}{2}\int_{R_z}(\delta_{\zeta}(X))^2d\zeta
\biggr\},\;\;z\in\mathbb{T}.
\end{equation}
Note that
\[
 V_z=\exp\biggl\{-\int_{R_z}\delta_{\zeta}(X)dW_{\zeta}^Y+\frac{1}{2}\int_{R_z}(\delta_{\zeta}(X))^2d\zeta
\biggr\},\;\;z\in\mathbb{T},
\]
thus, 
\begin{equation}\label{V2}
 V_z=\exp\biggl\{-\int_{R_z}\delta_{\zeta}(X)d\big(\int_{R_{\zeta}} K_{\alpha,\beta}^{-1}(\zeta;\zeta') dY_{\zeta'}
\big)
+\frac{1}{2}\int_{R_z}(\delta_{\zeta}(X))^2d\zeta
\biggr\},\;\;z\in\mathbb{T}.
\end{equation}
\begin{lemma} \label{mylemma}
Let $V=(V_z,z\in\mathbb{T})$ be defined by~\eqref{V2} (or, equivalently, by~\eqref{V1}). Then 
$\mathbb{E}\big(V_{(T_1,T_2)}\big)=1$.
\end{lemma}
{\bf Proof}: Since $B^{\alpha,\beta}$ and $X$ are independent, then $W^B$ and $X$ are independent, which implies that 
one can define a standard Wiener sheet $W^B$ on a complete probability space $(\Omega_2,\mathcal{F}_2,P_2)$, define 
$X$ on a complete probability space $(\Omega_1,\mathcal{F}_1,P_1)$ and then consider the processes on a product 
probability space $(\Omega_1\times\Omega_2,\mathcal{F}_1\times\mathcal{F}_2,P_1\times P_2)$, with 
$W^B(\omega)=W^B(\omega_2)$ and $X(\omega)=X(\omega_1)$ for all 
$\omega=(\omega_1,\omega_2)\in\Omega_1\times\Omega_2$. Let 
$Z(\omega)=Z(\omega_1,\omega_2)=
\int_{\mathbb{T}}\delta_{\zeta}(X(\omega_1))dW_{\zeta}^B(\omega_2)$. Then, upon taking into account Corollary~\ref{Cor1}, it follows that for almost all (fixed) $\omega_1$, 
$Z(\omega_1,\cdot)$ is a Gaussian random variable with mean 0 and variance $\int_{\mathbb{T}}
\big[\delta_{\zeta}(X(\omega_1))\big]^2d\zeta$, and, thus, for almost all $\omega_1\in\Omega_1$, 
\[\mathbb{E}_{P_2}(V_{(T_1,T_2)}(\omega_1,\cdot))=\mathbb{E}_{P_2}\!\left[\exp\bigg\{-Z(\omega_1,\cdot)-\frac{1}{2}\int_{\mathbb{T}}
\big[\delta_{\zeta}(X(\omega_1))\big]^2d\zeta\bigg\}\right]=1,
 \]
which implies that  
\[
 \mathbb{E}(V_{(T_1,T_2)})=\int_{\Omega_1\times\Omega_2}
V_{(T_1,T_2)}(\omega_1,\omega_2)(P_1\times P_2)(d\omega_1,d\omega_2)=
\int_{\Omega_1} 1\,P_1(d\omega_1)=1.\;\;\Box
\]
\begin{theorem} \label{bayestheorem}
Consider observation model~\eqref{observation}, where the signal 
$X=(X_z,z\in\mathbb{T})$ has almost surely H\"older-continuous sample paths of order $(\lambda_0,\lambda_0)$, where 
$\frac{\max(\alpha,\beta)-\frac{1}{2}}{\lambda}<\lambda_0<\frac{1}{2}$,  
and 
$g$ satisfies conditions ($\mathrm{A}_1$)-($\mathrm{A}_2$).
 Let $\tilde{P}$ be a new probability measure on $(\Omega,\mathcal{F})$ given by:
\begin{equation}
 \frac{d\tilde{P}}{dP}=V_{(T_1,T_2)}\;\;a.s.(P)
\end{equation}
Then $\tilde{P}$ is equivalent to $P$ and, under $\tilde{P}$, \eqref{observation} holds a.s., $Y$ is a standard fBs with Hurst indices $(\alpha,\beta)$, 
$X$ has the same law as under $P$, and processes $X$ and $Y$ are 
independent under $\tilde{P}$. Moreover, the following ``spatial-fractional'' version of 
the Bayes' formula holds: For any $F\in C_b(\mathbb{R})$, 
\begin{equation}
 \mathbb{E}\big(F(X_z)\vert\mathcal{F}_z^Y\big)=
\frac{\tilde{\mathbb{E}}\big[F(X_z)V^{-1}_{(T_1,T_2)}\vert\mathcal{F}_z^Y\big]}{\tilde{\mathbb{E}}
\big[V^{-1}_{(T_1,T_2)}\vert\mathcal{F}_z^Y\big]}=
\frac{\tilde{\mathbb{E}}\big[F(X_z)V^{-1}_{z}\vert\mathcal{F}_z^Y\big]}{\tilde{\mathbb{E}}
\big[V^{-1}_{z}\vert\mathcal{F}_z^Y\big]} \;\;a.s.
\end{equation}
where $\tilde{\mathbb{E}}$ denotes the mathematical expectation under $\tilde{P}$, 
$\mathcal{F}_z^Y$ denotes the filtration generated by the observation process in the 
rectangle $R_z=[0,z_1]\times[0,z_2]\subset\mathbb{T}$, i.e. 
\[
 \mathcal{F}_z^Y:=\sigma(Y_{\zeta}:\underline{0}\prec \zeta\prec z),\;\;z\in\mathbb{T},
\]
and $V=(V_z,z\in\mathbb{T})$ is defined by~\eqref{V2} (in terms of~\eqref{delta}).
\end{theorem}
{\bf{Proof}:} The first part of the theorem follows at once from the multiparameter Girsanov-type theorem for the standard 
Wiener sheet (see e.g.~Theorem~1 in~\cite{bib:dozzi}, p.~89), Proposition~\ref{EKtheorem} and Lemma~\ref{mylemma}. 
To prove Bayes' formula, arguments similar to those constructed in~\cite{bib:korezlioglu2} for the two-parameter 
Wiener sheet observation noise can be used. The proof, to a large extent, follows the lines of standard arguments 
used in the one-parameter (fractional noise) case, thus, we will omit the details here. $\Box$ 
\subsection{Evolution equation for the optimal nonlinear filter along an arbitrary 
increasing 1-dim curve}\label{curvesection}
Here we present a stochastic evolution equation satisfied by the unnormalized optimal 
filter when its dynamics is tracked along an arbitrary monotone non-decreasing 
1-dim continuous curve $\Delta$ connecting the origin to the point $T=(T_1,T_2)$. 
By a monotone non-decreasing path we mean that $\Delta$ is nondecreasing (in the sense 
of partial ordering in the plane) in both $z_1$ and $z_2$ directions. 
For each $z\in \mathbb{T}\equiv[0,T_1]\times[0,T_2]$, let $z_{\Delta}$ be the ``smallest" point on $\Delta$ which is larger 
than or equal to $z$ with respect to the partial ordering $\succ$.
The path $\Delta$ divides domain $\mathbb{T}$ into two regions; the region below $\Delta$, which is denoted by $D_{1}^{\Delta}$, and the region above $\Delta$, denoted by $D_{2}^{\Delta}$. \\ Namely,  
$D_{1}^{\Delta}=\big\{\zeta\in \mathbb{T}: 
\zeta \odot \zeta_{\Delta}=\zeta_{\Delta}\big\} $ and 
$D_{2}^{\Delta}=\big\{\zeta\in\mathbb{T}: \zeta_{\Delta} \odot 
\zeta=\zeta_{\Delta}\big\} $, where we use the notation $a\odot b:=(a_1,b_2)$ for arbitrary 
$a=(a_1,a_2),b=(b_1,b_2)\in\mathbb{R}^2_+$ (as in Section~\ref{multiparametermartingale}).
\begin{definition} Let $(\mathcal{F}_z,z\in\mathbb{T})$ be a filtration satisfying conditions (F1)-(F4) of  Section~\ref{multiparametermartingale}. Suppose $\Delta$ is a monotone nondecreasing continuous 1-dim curve connecting 
the origin to point $T=(T_1,T_2)\in\mathbb{R}^2_+$. Then \\
i) A process $\phi=(\phi_z,z\in\mathbb{T})$ is called $\Delta$-adapted if 
$\phi_z$ is $\mathcal{F}_{z_{\Delta}}$-measurable for all $z\in\mathbb{T}$. \\
ii) A process $X=(X_z,z\in\mathbb{T})$ is called a $\Delta$-martingale if $X$ is 
$\Delta$-adapted and 
\[
 \mathbb{E}\left.\big[X\big(z,z'\big]\,\right\vert\,\mathcal{F}_{z_{\Delta}}\big]=0\;\;
{for\;\;all\;\;} \underline{0}\prec z\prec z'\prec T.
\]
\end{definition}
\begin{definition} Let $\mathcal{H}_{\Delta}$ be the space of processes  
 $\phi=(\phi_z,z\in\mathbb{T})$ satisfying the following conditions:\\
{\rm(a)} $\phi$ is a bimeasurable function of $(\omega,z)$;\\
{\rm(b)} $\int_{\mathbb{T}}\mathbb{E}\phi_z^2dz<\infty$;\\
{\rm($\mathrm{c}_{\Delta}$)} $\phi$ is $\Delta$-adapted.  
\end{definition}
For $\phi\in\mathcal{H}_{\Delta}$, define processes $\phi_{i}^{\Delta}
=(\phi_{iz}^{\Delta},z\in\mathbb{T})\in\mathcal{H}_i$, $i=1,2$ 
(see Definition~\ref{H0} for definitions of $\mathcal{H}_1,\mathcal{H}_2$), by: 
\[
 \phi_{1z}^{\Delta}=\left\{
\begin{array}{ll}
 \phi_z, & \mathrm{if}\; z\in D_1^{\Delta},\\
0, & \mathrm{otherwise};
\end{array}\right.\;\;\mathrm{and}\;\;
\phi_{2z}^{\Delta}=\left\{
\begin{array}{ll}
 \phi_z, & \mathrm{if}\; z\in D_2^{\Delta},\\
0, & \mathrm{otherwise}.
\end{array}\right.
\]
Then $\phi_z=\phi_{1z}^{\Delta}+\phi_{2z}^{\Delta}$ for almost all $z\in\mathbb{T}$ and 
one can construct stochastic integral $\int_{\mathbb{T}}\phi_z dW_z=(\phi\circ W)^{\Delta}_T$ 
for $\phi\in\mathcal{H}_{\Delta}$ and 
show the 
following properties for the resulting integral (see~\cite{bib:zakai_curve} for details):
\begin{proposition} \label{curve_proposition}
Let $\Delta$ be a monotone nondecreasing 1-dim continuous curve connecting the origin to the 
 final point $T$. Let $\phi\in\mathcal{H}_{\Delta}$ and define the stochastic integral of $\phi$ with respect 
to a standard Wiener sheet $(W_z,\mathcal{F}_z, z\in\mathbb{T})$ by 
\begin{equation}\label{integral_Delta}
 (\phi\circ W)_z^{\Delta}=(\phi_1^{\Delta}\circ W)_z+(\phi_2^{\Delta}\circ W)_z,\;\;z\in\mathbb{T},
\end{equation}
where the two stochastic integrals on the right-hand side of~\eqref{integral_Delta} were discussed earlier 
in Section~\ref{multiparametermartingale}. Then the integral has the following properties:\\
i) $(\phi\circ W)^{\Delta}$ is a $\Delta$-martingale;\\
ii) $(\phi\circ W)^{\Delta}$ is a one-parameter martingale on the path $\Delta$;\\
iii) If $\Delta$ and $\Delta'$ are two monotone nondecreasing paths connecting the origin to $T$ and 
both passing through a point 
$z_0\in\mathbb{T}$, and $\phi$ is both $\Delta$ and $\Delta'$-adapted, then 
$(\phi\circ W)_{z_0}^{\Delta}=(\phi\circ W)_{z_0}^{\Delta'}$. 
\end{proposition}
 
Suppose our signal $X_{z}$ is a two-parameter semimartingale in the plane of the form: 
\begin{eqnarray}
X_{z}&=&X_{0}+\int_{R_{z}}\phi_{\zeta}dW_{\zeta}+\int_{R_{z}}\theta_{\zeta}d\zeta
+\int_{R_{z}\times R_{z}}\psi_{\zeta,\zeta'}dW_{\zeta}dW_{\zeta'}\\
&+&\int_{R_{z}\times R_{z}}f_{\zeta,\zeta'}d\zeta dW_{\zeta'}
+\int_{R_{z}\times R_{z}}g_{\zeta, \zeta'}dW_{\zeta}d\zeta',\;\;z\in\mathbb{T},
\end{eqnarray}
where, as usual, $R_z=[0,z_1]\times[0,z_2]$ and $\phi\in\mathcal{H}_{0}$ and
$\psi,f,g\in\hat{\mathcal{H}}$, where spaces $\mathcal{H}_0,\hat{\mathcal{H}}$ are defined as in 
Definition~\ref{H0} and Definition~\ref{hatH}.  

Then, by~\cite{bib:zakai_curve}, for an arbitrary monotone nondecreasing continuous 
1-dim curve $\Delta$, connecting the origin to $T$, there exist 
$\eta_{\zeta}=\eta(\Delta,\zeta)$ and 
$\nu_{\zeta}=\nu(\Delta, \zeta)$ such that $\eta\in \mathcal{H}_{\Delta}$ and 
\begin{equation}\label{signalsemimartingale}
X_{z}=X_{0}+\int_{R_{z}}\eta(\Delta,\zeta)dW_{\zeta}+\int_{R_{z}}\nu(\Delta,\zeta)d\zeta, \;\;z\in\Delta.
\end{equation}
If $\theta$ is $\Delta$-adapted, then $\nu$ can be chosen $\Delta$-adapted. As such, 
$X$ is clearly a sample-continuous semimartingale on $\Delta$.

In the rest of Section~\ref{curvesection} we will therefore assume that the signal process $X$ is of the 
form~\eqref{signalsemimartingale}, where the standard Wiener sheet $W$ is independent of the observation 
random field $Y$. Let us consider the nonlinear filtering model~\eqref{observation} along with conditions 
$\mathrm{(A}_1\mathrm{),(A}_2)$ and recall the general framework of Section~\ref{sectionbayes}. 
\begin{theorem} Let $\Delta$ be an arbitrary monotone nondecreasing continuous 1-dim curve connecting 
the origin to the final point $T\in\mathbb{R}^2_+$.  
Let us assume that the observation model~\eqref{observation} holds, along with conditions 
$\mathrm{(A}_1\mathrm{),(A}_2)$, and suppose that the signal $X$ is 
a two-parameter semimartingale in the plane, which is written in 
the form~{\eqref{signalsemimartingale}}, where $\eta\in\mathcal{H}_{\Delta}$ and $\nu$ is $\Delta$-adapted, and 
whose trajectories are 
H\"older-continuous of order 
$(\lambda_0,\lambda_0)$, where 
$\lambda_0>\frac{\max\{\alpha,\beta\}-\frac{1}{2}}{\lambda}$.  
For $F\in C_b^2(\mathbb{R})$, consider the unnormalized optimal filter 
\begin{equation}
\sigma_z(F):=\tilde{\mathbb{E}}\big[F(X_z)V_z^{-1}\vert\mathcal{F}_z^Y\big],\;\;z\in\mathbb{T},
\end{equation}
introduced in Theorem~{\ref{bayestheorem}}. Then the following stochastic evolution equation,  
governing the dynamics of the unnormalized optimal filter along the monotone increasing 
path $\Delta$, is satisfied:
\begin{equation}\label{zakai_curve}
 \sigma_z(F)=\sigma_{\underline{0}}(F)+\int_{R_z}
\sigma_{\zeta_{\Delta}}\big(\nu F'+\frac{1}{2}\eta^2 F''\big)d\zeta+
\int_{R_z}\sigma_{\zeta_{\Delta}}(F\delta)
d\biggl(\int_{R_{\zeta}}K_{\alpha,\beta}^{-1}(\zeta;\zeta')dY_{\zeta'}\biggr),\;\;z\in\Delta,
\end{equation}
where 
\begin{equation}\label{memoryeffect1}
 \sigma_{\zeta_{\Delta}}(F\delta):=\tilde{\mathbb{E}}\big[
F\left(X_{\zeta_{\Delta}}\right)\delta_{\zeta_{\Delta}}(X)V_{\zeta_{\Delta}}^{-1}\,\big\vert\,
\mathcal{F}_{\zeta_{\Delta}}^Y\big],
\end{equation}
and $(\delta_z(X),z\in\mathbb{T})$ is defined in~\eqref{delta}.
\end{theorem}
{\bf Proof:} 
Let us reparameterize $\Delta$ by $\{z(t);0\le t\le1\}$ so that the process 
$\{X_{z}, z\in \Delta \}$ can be rewritten as $\{X_{z(t)}, 0\le t \le1\}$. 
By Proposition~{\ref{curve_proposition}}, $X$ is a continuous one-parameter semimartingale 
on $\Delta$, thus, by It\^o's formula (for one-parameter case), for all $F\in C_b^2(\mathbb{R})$,
\[F(X_{z(t)})=F(X_{z(0)})+\int_{0}^{t}F'(X_{z(s)})dX_{z(s)}+\frac{1}{2}\int_{0}^{t}
F''(X_{z(s)})d\langle X, X\rangle _{z(s)}, \;\;t\in[0,1],\]
where
$\langle X,X\rangle_{z(t)}=\int_{R_{z(t)}}\eta_{\zeta}^{2}d\zeta$.
Note that one can re-express $F$ along $\Delta$ free of the earlier parametrization as follows: 
\[F(X_{z})=F(X_{\underline{0}})+\int_{R_{z}}F'(X_{\zeta_{\Delta}})dX_{\zeta}+
\frac{1}{2}\int_{R_{z}}F''(X_{\zeta_{\Delta}})\eta_{\zeta}^{2}d\zeta,   \quad  z\in \Delta.\]
Similarly, since
\[V_{z(t)}^{-1}=\exp\bigg\{\int_{0}^{t}\delta_{z(s)}(X)dW_{z(s)}^B+
\frac{1}{2}\int_{0}^{t}[\delta_{z(s)}(X)]^{2}dz(s)\bigg\},\;\;t\in[0,1],\]
then
\[
 V_{z(t)}^{-1}=1+\int_0^t V_{z(s)}^{-1}\delta_{z(s)}(X)dW^B_{z(s)}+
\int_0^t V_{z(s)}^{-1}[\delta_{z(s)}(X)]^2 dz(s),\;\;t\in[0,1],
\]
where the latter equation can also be rewritten free of parametrization as
\[V_{z}^{-1}=1+\int_{R_{z}}V_{\zeta_{\Delta}}^{-1}\delta_{\zeta_{\Delta}}(X)
dW_{\zeta}^Y,\;\;z\in \Delta.\]
Moreover,
\[V_{z(t)}^{-1}F(X_{z(t)})=F(X_{z(0)})+\int_{0}^{t}V_{z(s)}^{-1}\bigg(\nu_{z(s)}F'(X_{z(s)})+
\frac{1}{2}\eta_{z(s)}^{2}F''(X_{z(s)})\bigg)dz(s)\]
\[+\int_{0}^{t}V_{z(s)}^{-1}F'(X_{z(s)})\eta_{z(s)}dW_{z(s)}+
\int_{0}^{t}F(X_{z(s)})V_{z(s)}^{-1}\delta_{z(s)}(X)d{W}^Y_{z(s)},\;\;t\in[0,1].\]
Then, upon taking conditional expectations of both sides of the above one-parameter equation with respect 
to $\mathcal{F}_{z(t)}^Y=\mathcal{F}_{z(t)}^{W^Y}$ under $\tilde{P}$, one arrives at the  
following equation along the path $\Delta$:
\[\sigma_{z(t)}(F)=\sigma_{z(0)}(F)+\int_{0}^{t}\sigma_{z(s)}
(\nu F'+\frac{1}{2}\eta^{2}F'')dz(s)+\int_{0}^{t}\sigma_{z(s)}(F\delta)dW^Y_{z(s)},\;\;t\in[0,1],\]
where 
\[\sigma_{z(s)}(F\delta):=\tilde{E}[V_{z(s)}^{-1}F(X_{z(s)})\delta_{z(s)}(X)|\mathcal{F}_{z(s)}^Y].\]
The latter evolution along the 1-dimensional path $\Delta$ can be expressed free of parametrization as 
follows: 
 \[\sigma_{z}(F)=\sigma_{\underline{0}}(F)+\int_{R_{z}}\sigma_{\zeta_{\Delta}}
(\nu F'+\frac{1}{2}\eta^{2}F'')d\zeta
+\int_{R_{z}}\sigma_{\zeta_{\Delta}}(F\delta)dW^Y_{\zeta},\;\;z\in\Delta, \]
or, equivalently, 
\begin{equation}\label{multiparameterzakai}
\sigma_z(F)=\sigma_{\underline{0}}(F)+\int_{R_{z}}\sigma_{\zeta_{\Delta}}
(\nu F'+\frac{1}{2}\eta^{2}F'')d\zeta
+\int_{R_{z}}\sigma_{\zeta_{\Delta}}(F\delta)d\Big(\int_{R_{\zeta}}K_{\alpha,\beta}^{-1}
(\zeta;\zeta ')dY_{\zeta '}\Big),\;\;z\in\Delta,
\end{equation}
where the equations hold almost surely under $\tilde{P}$ and $P$. $\Box$
\begin{note} Let us observe that, in contrast to the case of filtering in the presence of a martingale observation 
 noise, the above stochastic evolution equation~\eqref{zakai_curve} cannot be interpreted  
as a measure-valued SPDE in 
view of the special meaning assigned to $\sigma(\cdot)$ in~\eqref{memoryeffect1}. The latter is necessary because   
$\delta_z(X)$ is not a function of $X_z$ but rather is a 
function of the entire ``history'' $(X_{\zeta},\underline{0}\prec \zeta\prec z)$.  
\end{note}

\subsection{Analogue of Duncan-Mortensen-Zakai equation for the optimal filter in 
the case of 2-parameter dynamics with fractional 
Brownian sheet noise}
The stochastic evolution equation~\eqref{zakai_curve} developed in Section~\ref{curvesection}, governing the dynamics of the unnormalized 
optimal filter, is two-dimensional in form, but clearly one-dimensional in spirit. Our objective in this section is to develop 
a ``fractional-spatial'' analogue of the  Duncan-Mortensen-Zakai equation for the unnormalized optimal filter which is 
inherently two-dimensional.  

Let $\mathfrak{a}:\mathbb{R}\rightarrow\mathbb{R}$ and $\mathfrak{b}:\mathbb{R}\rightarrow\mathbb{R}$ be measurable functions satisfying 
the following Lipshitz and growth conditions: there exists a finite constant $C>0$ such that for all $x,y\in\mathbb{R}$,
\[
 \vert \mathfrak{a}(x)-\mathfrak{a}(y)\vert+\vert\mathfrak{b}(x)-\mathfrak{b}(y)\vert\leq C\vert x-y\vert
\]
and 
\[
 \vert\mathfrak{a}(x)\vert+\vert\mathfrak{b}(x)\vert\leq C(1+\vert x\vert).
\]
Then there exists a unique strong solution to the following multiparameter SDE (see e.g.~\cite{bib:dozzi}): 
\[
X_z=X_{\underline{0}}+\int_{R_z} \mathfrak{a}(X_{\zeta})d\zeta+\int_{R_z}\mathfrak{b}(X_\zeta)dW_{\zeta},\;\;z\in\mathbb{T},
\]
where $W$ denotes a standard Wiener sheet. Moreover, the solution has H\"older-continuous sample path of order $(\lambda_1,\lambda_2)$ 
for all $\lambda_1,\lambda_2\in\big(0,\frac{1}{2}\big)$. 
\begin{theorem} \label{2d-evolution} In the framework of Section~\ref{sectionbayes}, assume that the observation model~\eqref{observation} holds, and that above Lipshitz and growth conditions on $\mathfrak{a}(\cdot)$ and $\mathfrak{b}(\cdot)$ are satisfied and  
conditions 
$\mathrm{(A}_1\mathrm{),(A}_2)$ are valid. 
Suppose that the signal $X$ is the unique strong solution of the following SDE:
\begin{equation}\label{diffusionsignal}
X_z=X_{\underline{0}}+\int_{R_z} \mathfrak{a}(X_{\zeta})d\zeta+\int_{R_z}\mathfrak{b}(X_\zeta)dW_{\zeta},\;\;z\in\mathbb{T},
\end{equation}
where $W$ is a standard Wiener sheet independent of the observation $Y$. Let $\sigma_z(F):=
\tilde{\mathbb{E}}\big[F(X_z)V_z^{-1}\,\vert\,\mathcal{F}_z^Y\big]$, i.e. $\sigma_z(F)$ is the unnormalized conditional expectation corresponding to the optimal filter. 
Then for all $F\in C_b^4(\mathbb{R})$, evolution of the unnormalized optimal filter has the following structure:
\[
 \sigma_{z}(F)= \sigma_{\underline{0}}(F) + \int_{R_z} \sigma_{\zeta}\big(\mathfrak{a}F'+\frac{1}{2}\mathfrak{b}^2F''\big)d\zeta+\int_{R_z}\sigma_{\zeta} 
\big(F\delta\big)d{W}_{\zeta}^Y
\]
\[
 +\iint_{R_z\times R_z}\sigma_{\zeta,\zeta'}\big(F;\delta\otimes\delta\big)
dW_{\zeta}^YdW_{\zeta'}^Y
\]
\[
 +\iint_{R_z\times R_z}\big[\sigma_{\zeta,\zeta'}(F';\mathfrak{a}\otimes\delta)+
\frac{1}{2}\sigma_{\zeta,\zeta'}(F'';\mathfrak{b}^2\otimes\delta)\big]d\zeta dW_{\zeta'}^Y
\]
\[
 +\iint_{R_z\times R_z}\big[\sigma_{\zeta,\zeta'}(F';\delta\otimes\mathfrak{a})+
\frac{1}{2}\sigma_{\zeta,\zeta'}(F'';\delta\otimes\mathfrak{b}^2)\big] dW_{\zeta}^Yd\zeta'
\]
\begin{equation}\label{evolution}
 +\!\iint_{R_z\times R_z}\!I(\zeta\curlywedge\zeta')\biggl[\sigma_{\zeta,\zeta'}\big(F'';\mathfrak{a}\otimes \mathfrak{a}\big)+
\frac{1}{2}\sigma_{\zeta,\zeta'}\big({F'''};\mathfrak{b}^2\otimes \mathfrak{a}+\mathfrak{a}\otimes 
\mathfrak{b}^2\big)
+\frac{1}{4}\sigma_{\zeta,\zeta'}\big(F^{(iv)};\mathfrak{b}^2\otimes \mathfrak{b}^2\big)
\biggr]d\zeta d\zeta',
\end{equation}
where $\otimes$ denotes the tensor product of functions, $\sigma_z(F\delta):=
\tilde{\mathbb{E}}\big[F(X_z)\delta_z(X)V_z^{-1}\,\vert\,\mathcal{F}_z^Y\big]$,\linebreak 
$\sigma_{z,z'}\big(F;\delta\otimes\delta\big):=
\tilde{\mathbb{E}}\big[F(X_{z\vee z'})\delta_z(X)\delta_{z'}(X)
V_{z\vee z'}^{-1}\,\vert\,\mathcal{F}_{z\vee z'}^Y\big]$, 
 and for arbitrary functions $f_1:\mathbb{R}\rightarrow\mathbb{R}$, $f_2:\mathbb{R}^2\rightarrow\mathbb{R}$, we put 
$\sigma_{z,z'}(f_1;f_2):=\tilde{\mathbb{E}}\big[f_1(X_{z\vee z'})f_2(X_z,X_{z'})
V_{z\vee z'}^{-1}\,\vert\,\mathcal{F}_{z\vee z'}^Y\big]$ for all $z,z'\in\mathbb{T}$. (In~\eqref{evolution}, $W_z^Y=\int_{R_z}
K_{\alpha,\beta}^{-1}(z;\zeta)dY_{\zeta}$ and $\delta$ is given by~\eqref{delta}, as before.)
\end{theorem}
\begin{note}
In Theorem~\ref{2d-evolution}, we could write $\sigma_z(F)=\sigma_{z,z}(F;1)$, where $1$ denotes function on 
$\mathbb{R}^2$ which is identically equal to one. 
\end{note}
{\bf Proof of Theorem~\ref{2d-evolution}}: First note that, under $\tilde{P}$, $Y$ is a fractional Brownian sheet with Hurst 
indices $(\alpha,\beta)$, 
while the corresponding field $W^Y$, given by 
$W^Y_z=\int_{R_z}K_{\alpha,\beta}^{-1}(z;\zeta)dY_{\zeta}$, is a standard Wiener sheet and the two random fields generate 
the same natural filtration, thus, the observation sigma-field $(\mathcal{F}_z^Y)_{0\prec z\prec T}$ has properties $(F1)$--$(F4)$ 
of Section~\ref{multiparametermartingale}. 
Similarly, $(\mathcal{F}_z^X)_{z\in\mathbb{T}}$ and $(\mathcal{F}_z^{X,Y})_{z\in\mathbb{T}}$ have properties 
$(F1)$--$(F4)$ under reference probability measure~$\tilde{P}$. 
Note also that the paths of $X=(X_z,z\in\mathbb{T})$ are almost surely H\"older-continuous of arbitrary 
order~$(\lambda_1,\lambda_2)$, where $\lambda_1,\lambda_2<\frac{1}{2}$, thus $\delta$ is well-defined and the conclusions of
Corollary~\ref{Cor1} hold. Next, by a version of the It\^o's formula for multiparameter semimartingales (see~\cite{bib:zakai2}), we 
obtain that for arbitrary $F\in C_b^4(\mathbb{R})$, 
\[
 F(X_z)=F(X_{\underline{0}})+\int_{R_z} F'(X_{\zeta})\big[
\mathfrak{a}(X_{\zeta})d\zeta+\mathfrak{b}(X_{\zeta})dW_{\zeta}\big]+\frac{1}{2}\int_{R_z}F''(X_{\zeta})\mathfrak{b}^2(X_{\zeta})
d\zeta
\]
\[
 +\iint_{R_z\times R_z}F''(X_{\zeta\vee\zeta'})\mathfrak{b}(X_{\zeta})\mathfrak{b}(X_{\zeta'})dW_{\zeta}dW_{\zeta'}
\]
\[
 +\iint_{R_z\times R_z}\big[F''(X_{\zeta\vee\zeta'})\mathfrak{b}(X_{\zeta})\mathfrak{a}(X_{\zeta'})+\frac{1}{2}
F'''(X_{\zeta\vee\zeta'})\mathfrak{b}(X_{\zeta})\mathfrak{b}^2(X_{\zeta'})\big]d\zeta dW_{\zeta'}
\]
\[
+\iint_{R_z\times R_z}\big[F''(X_{\zeta\vee\zeta'})\mathfrak{b}(X_{\zeta})\mathfrak{a}(X_{\zeta'})+\frac{1}{2}
F'''(X_{\zeta\vee\zeta'})\mathfrak{b}(X_{\zeta})\mathfrak{b}^2(X_{\zeta'})\big]dW_{\zeta}d\zeta'
\]
\[
 +\iint_{R_z\times R_z} I(\zeta\curlywedge\zeta')\biggl[
F''(X_{\zeta\vee\zeta'})\mathfrak{a}(X_{\zeta})\mathfrak{a}(X_{\zeta'})+
\frac{1}{2}F'''(X_{\zeta\vee\zeta'})\big(\mathfrak{a}(X_{\zeta})\mathfrak{b}^2(X_{\zeta'})+
\mathfrak{a}(X_{\zeta'})\mathfrak{b}^2(X_{\zeta})\big)
\]
\[
+\frac{1}{4}F^{(iv)}(X_{\zeta\vee\zeta'})\mathfrak{b}^2(X_{\zeta})\mathfrak{b}^2(X_{\zeta'})\biggr]
d\zeta d\zeta'.
\]
Similarly, under $\tilde{P}$, one shows that  
\[
 V_z^{-1}=1+\int_{R_z}V_{\zeta}^{-1}\delta_{\zeta}(X)dW_{\zeta}^Y+
\iint_{R_z\times R_z}V_{\zeta\vee\zeta'}^{-1}\,\delta_{\zeta}(X)\delta_{\zeta'}(X)dW_{\zeta}^YdW_{\zeta'}^Y\;\;
\mathrm{a.s.}
\]
Then the multiparameter version of the stochastic integration-by-parts formula (together with independence of $W$ and $W^Y$ under $\tilde{P}$) yields a corresponding equation for the product $F(X_z)V_z^{-1}$. Upon taking conditional expectation of both sides 
of the latter equation for $F(X_z)V_z^{-1}$ with respect to $\mathcal{F}^Y_z$ (where note that $\mathcal{F}^Y_z=\mathcal{F}_z^{W^Y}$)  
and using Lemma~{\ref{lemma-2d-evolution}}, which is proved below, one arrives at the following 
equation:
\[
\tilde{\mathbb{E}}\big(F(X_z)V_z^{-1}\vert\mathcal{F}_z^Y\big)=
\tilde{\mathbb{E}}\big(F(X_{\underline{0}})\vert\mathcal{F}_{\underline{0}}^Y\big)+
\int_{R_z}\tilde{\mathbb{E}}\biggl(\!\big[\mathfrak{a}(X_{\zeta})F'(X_{\zeta})+\frac{1}{2}
\mathfrak{b}^2(X_{\zeta})F''(X_{\zeta})\big]V_{\zeta}^{-1}\big\vert\mathcal{F}_z^Y\biggr)d\zeta
\]
\[
 +\int_{R_z}\tilde{\mathbb{E}}\big(F(X_{\zeta})\delta_{\zeta}(X)V_{\zeta}^{-1}\,\vert\,\mathcal{F}_{\zeta}^Y\big)dW_{\zeta}^Y\]
\[
+\iint_{R_z\times R_z}
\tilde{\mathbb{E}}\biggl(F(X_{\zeta\vee\zeta'})\delta_{\zeta}(X)\delta_{\zeta'}(X)V_{\zeta\vee\zeta'}^{-1}\,\biggr\vert\,
\mathcal{F}_{\zeta\vee\zeta'}^Y\biggr)dW^Y_{\zeta}dW^Y_{\zeta'}
\]
\[
 +\iint_{R_z\times R_z}
\tilde{\mathbb{E}}\biggl(\big[\mathfrak{a}(X_{\zeta})F'(X_{\zeta\vee\zeta'})+\frac{1}{2}\mathfrak{b}^2(X_{\zeta})
F''(X_{\zeta\vee\zeta'})\big]\delta_{\zeta'}V_{\zeta\vee\zeta'}^{-1}\bigg\vert\mathcal{F}_{\zeta\vee\zeta'}^Y\biggr)d\zeta dW^Y_{\zeta'}
\]
\[
 +\iint_{R_z\times R_z}
\tilde{\mathbb{E}}\biggl(\big[\mathfrak{a}(X_{\zeta'})F'(X_{\zeta\vee\zeta'})+\frac{1}{2}\mathfrak{b}^2(X_{\zeta'})
F''(X_{\zeta\vee\zeta'})\big]\delta_{\zeta}V_{\zeta\vee\zeta'}^{-1}\bigg\vert\mathcal{F}_{\zeta\vee\zeta'}^Y\biggr)
dW^Y_{\zeta} d\zeta'
\]
\[
+\iint_{R_z\times R_z}\!I(\zeta\curlywedge\zeta')
\tilde{\mathbb{E}}\biggl(\!\biggl[F''(X_{\zeta\vee\zeta'})\mathfrak{a}(X_{\zeta})\mathfrak{a}(X_{\zeta'})
+\frac{1}{4}F^{(iv)}(X_{\zeta\vee\zeta'})\mathfrak{b}^2(X_{\zeta})\mathfrak{b}^2(X_{\zeta'})
\]
\[
+\frac{1}{2}F'''(X_{\zeta\vee\zeta'})\big\{\mathfrak{b}^2(X_{\zeta'})\mathfrak{a}(X_{\zeta})+
\mathfrak{a}(X_{\zeta'})\mathfrak{b}^2(X_{\zeta})\big\}
\biggr]
V_{\zeta\vee\zeta'}^{-1}\,\biggr\vert\,
\mathcal{F}_{\zeta\vee\zeta'}^Y\biggr)d\zeta d\zeta'\;\;\mathrm{a.s.},
\]
thus, the required conclusion follows. $\Box$
\begin{lemma}\label{lemma-2d-evolution} Let $W$ and $W^Y$ be independent standard Wiener sheets on a probability space 
 $(\Omega,\mathcal{F}_T,\tilde{P})$ and $\mathcal{F}^{W,W^Y}_z:=\sigma(W_{\zeta}, W^Y_{\zeta'}:\underline{0}\prec \zeta\prec z, 
0\prec\zeta'\prec z)$, $z\in\mathbb{T}$. Also let $(\mathcal{F}_z^W)$ and $(\mathcal{F}_z^{W^Y})$ denote the natural filtrations generated 
by $W$ and $W^Y$, respectively. 
Consider 
a process $M$ (which is $\big(\mathcal{F}_z^{W, W^Y}\big)$-measurable), given by  
\[
 M_z:=\int_{R_z}\phi_{\zeta}dW_{\zeta}^Y+\iint_{R_z\times R_z}\psi_{\zeta,\zeta'} dW_{\zeta}^YdW_{\zeta'}^Y,
\]
where $\phi\in\mathcal{H}_0$ and $\psi\in\hat{\mathcal{H}}$, with $\mathcal{H}_0$ and $\hat{\mathcal{H}}$ being defined with respect 
to filtration $(\mathcal{F}_z^{W, W^Y})$. Then 

i) For any process $\psi\in\hat{\mathcal{H}}$, 
\[
 \tilde{\mathbb{E}}\left(\int_{R_z\times R_z}\psi_{\zeta,\zeta'}dW_{\zeta}dW_{\zeta'}\,\biggr\vert\,\mathcal{F}_z^{W^Y}\right)=0\;\;
\mathrm{a.s.}\;\;\tilde{P},
\]
\[
 \tilde{\mathbb{E}}\left(\int_{R_z\times R_z}\psi_{\zeta,\zeta'}dW_{\zeta}d\zeta'\,\biggr\vert\,\mathcal{F}_z^{W^Y}\right)=0\;\;
\mathrm{a.s.}\;\;\tilde{P},
\]
\[
 \tilde{\mathbb{E}}\left(\int_{R_z\times R_z}\psi_{\zeta,\zeta'}d\zeta dW_{\zeta'}\,\biggr\vert\,\mathcal{F}_z^{W^Y}\right)=0\;\;
\mathrm{a.s.}\;\;\tilde{P},
\]
\[
 \tilde{\mathbb{E}}\left(\int_{R_z\times R_z}\psi_{\zeta,\zeta'}dW_{\zeta}dW_{\zeta'}^Y\,\biggr\vert\,\mathcal{F}_z^{W^Y}\right)=0\;\;
\mathrm{a.s.}\;\;\tilde{P},
\]
\[
 \tilde{\mathbb{E}}\left(\int_{R_z\times R_z}\psi_{\zeta,\zeta'}dW_{\zeta}^YdW_{\zeta'}\,\biggr\vert\,\mathcal{F}_z^{W^Y}\right)=0\;\;
\mathrm{a.s.}\;\;\tilde{P}.
\]
ii) The following equation holds almost surely with respect to $\tilde{P}$:
\[
 \tilde{\mathbb{E}}\left(M_z\,\vert\,\mathcal{F}_z^{W^Y}\right)=
\int_{R_z}\tilde{\mathbb{E}}\left(\phi_{\zeta}\,\vert\,
\mathcal{F}_{\zeta}^{W^Y}\right)dW_{\zeta}^Y+\int_{R_z\times R_z}
\tilde{\mathbb{E}}\left(\psi_{\zeta,\zeta'}\,\vert\,
\mathcal{F}_{\zeta\vee\zeta'}^{W^Y}\right)dW_{\zeta}^Y dW_{\zeta'}^Y.
\]
Also, 
\[
 \tilde{\mathbb{E}}\left(\int_{R_z\times R_z}\psi_{\zeta,\zeta'}d\zeta dW_{\zeta'}^Y\,\biggr\vert\,\mathcal{F}_z^{W^Y}\right)=
\int_{R_z\times R_z}\tilde{\mathbb{E}}\big(
\psi_{\zeta,\zeta'}\vert\mathcal{F}_{\zeta\vee\zeta'}^{W^Y}\big)d\zeta dW_{\zeta'}^Y
\;\;
\mathrm{a.s.}\;\;\tilde{P},
\]
\[
 \tilde{\mathbb{E}}\left(\int_{R_z\times R_z}\psi_{\zeta,\zeta'}dW_{\zeta}^Yd\zeta'\,\biggr\vert\,\mathcal{F}_z^{W^Y}\right)=
\int_{R_z\times R_z}\tilde{\mathbb{E}}\big(
\psi_{\zeta,\zeta'}\vert\mathcal{F}_{\zeta\vee\zeta'}^{W^Y}\big)dW_{\zeta}^Y d\zeta'
\;\;
\mathrm{a.s.}\;\;\tilde{P}.
\]
\end{lemma}
{\bf Proof}: Let us start by showing that the first equality in~(i) holds, i.e. that 
\[\tilde{\mathbb{E}}\bigg(\int_{R_z\times R_z}\psi_{\zeta,\zeta'}dW_{\zeta}dW_{\zeta'} 
\vert\mathcal{F}_{z}^{W^Y}\bigg)=0
\] almost surely 
under $\tilde{P}$.  By independence of $W$ and $W^Y$, we may assume that $W$ is a standard Wiener sheet on a filtered complete probability space $(\Omega^X,\breve{\mathcal{F}}^W_T,(\breve{\mathcal{F}}_z^W)_{z\in\mathbb{T}}, \tilde{P}^X)$, whereas $W^Y$ is a standard Wiener sheet on another  filtered complete probability space $(\Omega^Y,\breve{\mathcal{F}}^{W^Y}_T,(\breve{\mathcal{F}}_z^{W^Y})_{z\in\mathbb{T}}, \tilde{P}^Y)$, 
where $(\breve{\mathcal{F}}_z^W)$ and $(\breve{\mathcal{F}}_z^{W^Y})$ are, respectively, natural filtrations generated by processes $W$ and $W^Y$ (in $\Omega^X$ and $\Omega^Y$, respectively), 
and  
$(\Omega,\mathcal{F}_T,\tilde{P})=(\Omega^X\times\Omega^Y,\breve{\mathcal{F}}^{W}_T\times\breve{\mathcal{F}}^{W^Y}_T,P^X\times P^{Y})$, 
i.e. the product probability space. Then $W$ and $W^Y$ are defined on $(\Omega,\mathcal{F}_T,\tilde{P})$ by 
$W_z(\omega)=W_z(\omega_1)$ and $W_z^Y(\omega)=W_z^Y(\omega_2)$ for all $\omega=(\omega_1,\omega_2)\in\Omega$. Then, clearly, 
$\mathcal{F}_z^W=\breve{\mathcal{F}}_z^W\times\{\emptyset,\Omega^Y\}$ and $\mathcal{F}_z^{W^Y}=
\{\emptyset,\Omega^X\}\times\breve{\mathcal{F}}_z^{W^Y}$. 

Next let us fix some $n\in\mathbb{N}$ and consider a partition of rectangle $R_T=(\underline{0},T]$ (where $T=(T_1,T_2)$) into 
rectangles $\Delta_{i,j}:=\big(z_{(i,j)},z_{(i+1,j+1)}\big]$, where $z_{(i,j)}=(2^{-n}iT_1, 2^{-n}jT_2)$. Let $\mathcal{S}$ be the 
class of processes $\psi$ of the form:
\begin{equation}\label{simple_lemma}
 \psi(\zeta,\zeta')=\sum_{i,j,k,\ell=0}^{2^n-1}\alpha_{ijk\ell}1_{\Delta_{i,j}}(\zeta)1_{\Delta_{k,\ell}}(\zeta'),
\end{equation}
where $\alpha_{ijk\ell}$ is $\mathcal{F}_{z_{(i,j)}\vee z_{(k,\ell)}}^{W,W^Y}$-measurable. By definition of the double integral,
\[
 \int_{R_z\times R_z}\psi(\zeta,\zeta')dW_{\zeta}dW_{\zeta'}:=\sum_{i,j,k,\ell=0}^{2^n-1}
\alpha_{ijk\ell}\,1_{\{i<k\}}1_{\{\ell<j\}}W(R_z\cap\Delta_{i,j})W(R_z\cap\Delta_{k,\ell}).
\]
Then, 
\begin{eqnarray*}
 & &\tilde{\mathbb{E}}\bigg(
\int_{R_z\times R_z}\psi(\zeta,\zeta')dW_{\zeta}dW_{\zeta'}
\vert\mathcal{F}_z^{W^Y}\bigg)\\
&=&\sum_{i,j,k,\ell=0}^{2^n-1}\tilde{\mathbb{E}}\biggr[\alpha_{ijk\ell}1_{\{i<k\}\cap\{\ell<j\}}
W(R_z\cap\Delta_{i,j})W(R_z\cap\Delta_{k,\ell})\,\big\vert\,\{\emptyset,\Omega^X\}\times\breve{\mathcal{F}}_z^{W^Y}\biggr].
\end{eqnarray*}
Note that $\forall Q\in\breve{\mathcal{F}}_z^{W^Y}$, we have 
\[
 \int_{\Omega^X\times Q}\alpha_{ijk\ell}(\omega_1,\omega_2)1_{\{i<k\}\cap\{\ell<j\}}
W(R_z\cap\Delta_{i,j})(\omega_1)W(R_z\cap\Delta_{k,\ell})(\omega_1)d\tilde{P}(\omega_1,\omega_2)
\]
\[
 =\!\int_Q\!\left[
\int_{\Omega^X}\alpha_{ijk\ell}(\omega_1,\omega_2)1_{\{i<k\}\cap\{\ell<j\}}
W(R_z\cap\Delta_{i,j})(\omega_1)W(R_z\cap\Delta_{k,\ell})(\omega_1)d\tilde{P}^X(\omega_1)\right]\!d\tilde{P}^Y(\omega_2)
\]
$=0$, since for fixed $\omega_2\in\Omega_2$ and for all $i<k$ and $\ell<j$, random variables 
$\alpha_{ijk\ell}(\cdot,\omega_2)$, $W(R_z\cap\Delta_{k,\ell})$ and $W(R_z\cap\Delta_{i,j})$ are mutually independent. Thus, 
\[
 \tilde{\mathbb{E}}\biggr[\alpha_{ijk\ell}1_{\{i<k\}\cap\{\ell<j\}}
W(R_z\cap\Delta_{i,j})W(R_z\cap\Delta_{k,\ell})\,\big\vert\,\{\emptyset,\Omega^X\}\times\breve{\mathcal{F}}_z^{W^Y}\biggr]=0\;\;
\mathrm{a.s.}(\tilde{P}),
\]
implying that for any simple process $\psi\in \mathcal{S}$, 
\[
 \tilde{\mathbb{E}}\bigg(
\int_{R_z\times R_z}\psi(\zeta,\zeta')dW_{\zeta}dW_{\zeta'}
\vert\mathcal{F}_z^{W^Y}\bigg)=0\;\;\mathrm{a.s.} \;(\tilde{P}).
\]
Since $\mathcal{S}$ is dense in $\hat{\mathcal{H}}$, the required equality follows for arbitrary $\psi\in\hat{\mathcal{H}}$ by 
taking appropriate limits.  
Similar arguments show that the remaining equalities in (i) are also valid. 

To prove (ii), let us show that $\forall\psi\in\hat{\mathcal{H}}$, 
\begin{equation}
 \tilde{\mathbb{E}}\bigg(
\int_{R_z\times R_z}\psi(\zeta,\zeta')dW_{\zeta}^YdW_{\zeta'}^Y
\vert\mathcal{F}_z^{W^Y}\bigg)=\int_{R_z\times R_z}\tilde{\mathbb{E}}\big[\psi(\zeta,\zeta')\vert\mathcal{F}^{W^Y}_{\zeta\vee\zeta'}\big]dW_{\zeta}^YdW_{\zeta'}^Y\;\;\mathrm{a.s.}
\label{lemma(ii)}
\end{equation}
First, consider $\psi\in\mathcal{S}$ of the form~{\eqref{simple_lemma}}. Then
\begin{eqnarray}
&& \tilde{\mathbb{E}}\bigg(
\int_{R_z\times R_z}\psi(\zeta,\zeta')dW_{\zeta}^YdW_{\zeta'}^Y
\vert\mathcal{F}_z^{W^Y}\bigg)\nonumber\\
&=&\sum_{i,j,k,\ell=0}^{2^n-1}\tilde{\mathbb{E}}\big[\alpha_{ijk\ell}\vert\mathcal{F}_z^{W^Y}\big]
1_{\{i<k\}\cap\{\ell<j\}}W^Y(R_z\cap\Delta_{i,j})W^Y(R_z\cap\Delta_{k\ell}) \label{convenience1},
\end{eqnarray}
where note that 
\[1_{\{i<k\}\cap\{\ell<j\}}W^Y(R_z\cap\Delta_{i,j})W^Y(R_z\cap\Delta_{k\ell})=0 \;\;\mathrm{unless}\;\;z_{(k,j)}=
\big(z_{(i,j)}\vee z_{(k,\ell)}\big)\prec\prec z.
 \]
Since $\alpha_{ijk\ell}$ is $\mathcal{F}^{W,W^Y}_{z_{(k,j)}}$-measurable and $z_{(k,j)}\prec\prec z$, then the conditional 
expectation in the right-hand side of~\eqref{convenience1} satisfies equation
\[
\tilde{\mathbb{E}}\big[\alpha_{ijk\ell}\vert\mathcal{F}_z^{W^Y}\big]=
\tilde{\mathbb{E}}\big[\alpha_{ijk\ell}\vert\mathcal{F}_{z_{(k,j)}}^{W^Y}\big]\;\;\mathrm{a.s.},
\]
by independence of $W$ and $W^Y$ and since Wiener sheets generate independently scattered measures.  
Thus, 
 \begin{eqnarray}
&& \tilde{\mathbb{E}}\bigg(
\int_{R_z\times R_z}\psi(\zeta,\zeta')dW_{\zeta}^YdW_{\zeta'}^Y
\vert\mathcal{F}_z^{W^Y}\bigg)\nonumber\\
&=&\sum_{i,j,k,\ell=0}^{2^n-1}\tilde{\mathbb{E}}\big[\alpha_{ijk\ell}\vert\mathcal{F}_{z_{(k,j)}}^{W^Y}\big]
1_{\{i<k\}\cap\{\ell<j\}}W^Y(R_z\cap\Delta_{i,j})W^Y(R_z\cap\Delta_{k\ell})\;\;\mathrm{a.s.} \label{convenience2}
\end{eqnarray}
On the other hand, 
\begin{eqnarray}
&& \int_{R_z\times R_z}\tilde{\mathbb{E}}\big[\psi(\zeta,\zeta')\vert\mathcal{F}^{W^Y}_{\zeta\vee\zeta'}\big]dW_{\zeta}^YdW_{\zeta'}^Y
\nonumber\\
&=&\int_{R_z\times R_z}\sum_{i,j,k,\ell=0}^{2^n-1}\tilde{\mathbb{E}}\big[\alpha_{ijk\ell}\vert\mathcal{F}^{W^Y}_{\zeta\vee\zeta'}\big]
1_{\Delta_{i,j}}(\zeta)1_{\Delta_{k,\ell}}(\zeta')dW_{\zeta}^Y dW_{\zeta'}^Y\nonumber\\
&=&\int_{R_z\times R_z}\sum_{i,j,k,\ell=0}^{2^n-1}1_{\{i<k\}\cap\{\ell<j\}}
\tilde{\mathbb{E}}\big[\alpha_{ijk\ell}\vert\mathcal{F}^{W^Y}_{\zeta\vee\zeta'}\big]
1_{\Delta_{i,j}}(\zeta)1_{\Delta_{k,\ell}}(\zeta')dW_{\zeta}^Y dW_{\zeta'}^Y, \label{convenience3}
\end{eqnarray}
where the last equality holds by definition of the double integral. 
Note that
\[
 1_{\{i<k\}\cap\{\ell<j\}}1_{\Delta_{i,j}}(\zeta)1_{\Delta_{k,\ell}}(\zeta')\neq 0\;\;\mathrm{implies\;\;that}\;\;
\zeta\vee\zeta'\in\Delta_{z_{(k,j)}},
\]
which, in turn, implies that $\tilde{\mathbb{E}}\big[\alpha_{ijk\ell}\vert\mathcal{F}_{\zeta\vee\zeta'}^{W^Y}\big]$ (in the right-hand side 
of~\eqref{convenience3}) equals almost surely to 
$\tilde{\mathbb{E}}\big[\alpha_{ijk\ell}\vert\mathcal{F}_{z_{(k,j)}}^{W^Y}\big]$, since $\alpha_{ijk\ell}$ is 
$\mathcal{F}_{z_{(k,j)}}^{W,W^Y}$-measurable. Thus, from~\eqref{convenience3} by definition of the double integral, 
\begin{eqnarray}
&& \int_{R_z\times R_z}\tilde{\mathbb{E}}\big[\psi(\zeta,\zeta')\vert\mathcal{F}^{W^Y}_{\zeta\vee\zeta'}\big]dW_{\zeta}^YdW_{\zeta'}^Y
\nonumber\\
&=&\sum_{i,j,k,\ell=0}^{2^n-1}\tilde{\mathbb{E}}\big[\alpha_{ijk\ell}\vert\mathcal{F}_{z_{(k,j)}}^{W^Y}\big]
1_{\{i<k\}\cap\{\ell<j\}}W^Y(R_z\cap\Delta_{i,j})W^Y(R_z\cap\Delta_{k\ell})\;\;\mathrm{a.s.} \label{convenience4}
\end{eqnarray}
From~\eqref{convenience2} and~\eqref{convenience4}, it follows that~\eqref{lemma(ii)} holds for all $\psi\in\mathcal{S}$. Since 
$\mathcal{S}$ is dense in $\hat{\mathcal{H}}$, it follows that~\eqref{lemma(ii)} holds for all $\psi\in\hat{\mathcal{H}}$ by 
taking appropriate limits. Similarly one establishes that $\forall\phi\in\mathcal{H}_0$,  
\[
\tilde{\mathbb{E}}\big(
\int_{R_z}\phi_{\zeta}dW_{\zeta}^Y\big\vert\mathcal{F}_z^{W^Y}\big)=\int_{R_z}\tilde{\mathbb{E}}(\phi_\zeta\vert
\mathcal{F}_{\zeta}^{W^Y})dW_{\zeta}^Y\;\;\mathrm{a.s.},
\]
thus, the first statement in (ii) is proved. The remaining two statements in (ii) can be established by analogous arguments.
$\Box$

\section{Conclusions}
In this paper the problem of spatial nonlinear filtering of a multiparameter semimartingale random 
field, with estimation based on an observation random field perturbed by a long-memory fractional noise, has been 
considered. Two types of stochastic evolution equations, governing the dynamics of the unnormalized optimal filter in the 
2-dimensional plane, has been derived. One equation follows the dynamics of the optimal filter along an arbitrary 
non-decreasing (in the sense of partial ordering) one-dimensional curve, while the other describes behavior of 
the optimal filter in terms of ``truly'' 2-dimensional dynamics. In view of long-memory in the observation noise, 
neither equation can be viewed as measure-valued SPDE and their interpretation is not trivial. However natural questions 
regarding uniqueness and robustness of the solutions to the evolution equations, as well as construction of suboptimal filters, can be addressed and the authors plan to do so  
in the forthcoming work.  

Despite numerous important practical applications of spatial nonlinear filtering (in connection with ``denoising'' and filtering of images and 
video-streams in physical, biological and atmospheric sciences, for example), there currently appears to be very little mathematical 
literature on the subject. In particular, the results presented in this paper represent the first mathematical 
results pertaining to spatial nonlinear filtering of random fields 
in the presence of long-memory (fractional) spatial observation noise.

\end{document}